\documentclass[lettersize,journal]{IEEEtran}
\usepackage{amsmath,amsfonts}
\usepackage{amssymb,amsthm}
\usepackage{algorithm}
\usepackage{algorithmicx,algorithm}
\usepackage{algpseudocode}
\usepackage{listings}
\usepackage{array}
\usepackage{textcomp}
\usepackage{stfloats}
\usepackage{url}
\usepackage{verbatim}
\usepackage{graphicx}
\usepackage{multirow}
\usepackage{subfigure}
\usepackage{epsfig}
\usepackage{epstopdf}
\usepackage{booktabs}
\usepackage[titles,subfigure]{tocloft}
\usepackage{float}
\usepackage{caption}
\usepackage{cite}
\usepackage{hyperref}
\usepackage{cleveref}
\Crefname{figure}{Fig.}{Figs.}
\usepackage{xcolor}
\definecolor{myBlue}{RGB}{0, 0, 255}
\hypersetup{
	colorlinks=true,
	linkcolor=myBlue,
	citecolor=myBlue,
	urlcolor=myBlue
}
\usepackage{orcidlink}
\begin{document}
	\theoremstyle{thmstylethree}%
	\newtheorem{definition}{Definition}%

\theoremstyle{thmstyleone}%
\newtheorem{theorem}{Theorem}

\title{Nonlinear kernel-free quadratic hyper-surface support vector machine with 0-1 loss function}
	\author{Mingyang Wu, Zhixia Yang*, Junyou Ye
\thanks{Mingyang Wu, Zhixia Yang and Junyou Ye are with the College of Mathematics and System Science, Xinjiang University, Urumqi 830046, China, (email:wmyang068@163.com; yangzhx@xju.edu.cn; yejymath@163.com).}
\thanks{*Corresponding author.}
}
\maketitle

\begin{abstract}
	For the binary classification problem, a novel nonlinear kernel-free quadratic hyper-surface support vector machine with 0-1 loss function (QSSVM$_{0/1}$) is proposed. Specifically, the task of QSSVM$_{0/1}$ is to seek a quadratic separating hyper-surface to divide the samples into two categories. And it has better interpretability than the methods using kernel functions, since each feature of the sample acts both independently and synergistically. By introducing the 0-1 loss function to construct the optimization model makes the model obtain strong sample sparsity. The proximal stationary point of the optimization problem is defined by the proximal operator of the 0-1 loss function, which figures out the problem of non-convex discontinuity of the optimization problem due to the 0-1 loss function.  
A new iterative algorithm based on the alternating direction method of multipliers (ADMM) framework is designed to solve the optimization problem, which relates to the working set defined by support vectors. The computational complexity and convergence of the algorithm are discussed. 
Numerical experiments on 4 artificial datasets and 14 benchmark datasets demonstrate that our QSSVM$_{0/1}$ achieves higher classification accuracy, fewer support vectors and less CPU time cost than other state-of-the-art methods. 
\end{abstract}


\section{Introduction}
\IEEEPARstart{S}{upport} vector machine (SVM) is one of the most popular and crucial machine learning methods based on the statistical learning theory and structural risk minimization principle. It was initially proposed by \cite{cortes1995support}. Since then, SVMs have been extensively studied and successfully applied in various domains, including text categorization \cite{joachims2005text}, pattern recognition \cite{10.1023/A:1009715923555}, image processing \cite{10.5555/1853979.1853985}, and biomedical \cite{2005Genetic}, etc.

For a given binary classification problem, SVM \cite{cortes1995support} aims to obtain a hyperplane to separate samples into two classes by solving a convex quadratic programming problem. 
As a key element in machine learning, the loss function is a non-negative real valued function used to constrain and guide model optimization. To enhance the performance of the model, various loss functions have been used in SVM-type methods. These methods include the least squares support vector machine \cite{1999Least, 1963An} (LSSVM), SVM with pinball loss \cite{huang2013support, MR3947658} (PSVM), SVM with ramp soft-margin loss \cite{collobert2006trading, Brooks2011Support, 2008Training} (RSVM), and SVM with non-convex robust and smooth soft-margin loss \cite{10.1162/NECO_a_00837} (RSHSVM), etc. However these existing soft-margin SVMs based on other loss functions can be viewed as substitutes for the 0-1 loss SVM. Although the 0-1 loss function is an ideal loss function, its application is limited, since it is discontinuous and non-convex. Recently, \cite{15921057020221001} proposed a SVM classifier via 0-1 soft-margin loss (SVM$_{0/1}$). Although the optimization problem is NP-hard problem \cite{1995Sparse}, research shows that it can be solved directly. And by using the 0-1 loss function it can achieve strong sparsity and robustness. This drives the development of the 0-1 loss function. However, SVM$_{0/1}$ can only obtain linear decision function and its application has limitations. Then, a nonlinear version \cite{liu2024nonlinear} was proposed, in which the decision function is generated by a kernel function. It can handle nonlinear data, but may lead to overfitting and lack of interpretability. Additionally, it is challenging to select the appropriate kernel functions and corresponding parameters. 

To overcome the above mentioned  limitations, \cite{MR2386593} developed a quadratic kernel-free nonlinear SVM based on the maximum margin idea, which paved the way for the development of the kernel-free trick. Subsequently, the soft-margin quadratic hyper-surface SVM \cite{MR3580868} (SQSSVM) was proposed to improve robustness and classification accuracy. The purpose of SQSSVM is to find a quadratic separating hyper-surface to classify samples into two classes. Following this, \cite{MR3411766} proposed the quadratic kernel-free least squares SVM. Furthermore, the kernel-free trick has also been used for semi-supervised \cite{11670570620160701} and unsupervised \cite{MR4011147} learning. \cite{sym13081378} proposed the kernel-free quadratic hyper-surface minimax probability machine (QSMPM) for a binary classification problem. Additional classification methods based on kernel-free trick can be found in the literature \cite{15737104020220501, MR4188956, 15856407920220915, MR4034610, MR4226386, TIAN201896,10.1007/s00500-017-2751-z, 10.1109/TFUZZ.2017.2752138}. Moreover, the kernel-free trick has also been applied to regression \cite{YE2022177, MR4083319} and clustering \cite{10.1007/s00500-016-2462-x} problems. Therefore, the kernel-free trick plays a significant role in machine learning.

In this paper, we develop a novel kernel-free quadratic hyper-surface method with 0-1 loss function for the binary classification problem, namely QSSVM$_{0/1}$. It is constructed based on the kernel-free trick and inherits the advantages of 0-1 loss function. Specifically, the main contributions of this paper are as follows.
\begin{itemize} \item QSSVM$_{0/1}$ is proposed by utilizing the kernel-free trick, it avoids the difficulty of selecting kernel functions and corresponding parameters. It has better interpretability than the methods using kernel functions, since each feature of the sample acts both independently and synergistically. 
	\item The 0-1 loss function is introduced to build the optimization problem to obtain better sample sparsity. By using the proximal operator of the 0-1 loss function to define the P-stationary point of the optimization problem, the non-convex discontinuity of the optimization problem due to the 0-1 loss function is solved.
	
	\item The support vectors are defined, and all support vectors fall on the support hyper-surfaces. Additionally, the relationship between the optimal solution to the  optimization problem and proximal stationary point is discussed. Under certain conditions, the proximal stationary point is the optimal solution.
	\item A working set is defined based on the information of the proximal operator of the 0-1 loss function. An iterative algorithm for solving the optimization problem of QSSVM$_{0/1}$ is given on the working set. In addition, the complexity and convergence of the algorithm are explored.
	\item  Numerical experiments are performed on 4 artificial datasets and 14 benchmark datasets. The experimental results show that our QSSVM$_{0/1}$ has higher classification accuracy and fewer support vectors than most methods. Moreover, the computational speed of our QSSVM$_{0/1}$ is fast. 
\end{itemize}

The rest of the paper is organized as follows. Section \ref{2} describes some symbols and reviews some related work for binary classification. Section \ref{3} presents our binary classification method, QSSVM$_{0/1}$. Numerical experiments on the artificial datasets and benchmark datasets are shown in section \ref{4}. Section \ref{5} concludes the paper.

	\section{Related work}\label{2}
In this section, we start with some notations and definitions, briefly introduce two classifying methods, SVM$_{0/1}$ and SQSSVM.

	\subsection{Symbol description}
Through this article, we use lower case letters to represent scalars, bold lower case letters to represent vectors, and bold upper case letters to represent matrices. $\mathbb{R}$ denotes the set of real numbers, $\mathbb{R}^{n}$ represents the $n$-dimensional real space, $\mathbb{S}^{n}$ represents the $n\times n$ dimensional real symmetric matrix.
For convenience, some operators are defined as follows.
\begin{definition}\label{defi1}
	The half-vectorization operator of any symmetric matrix $\boldsymbol{W}=[w_{ij}]_{n\times n}\in\mathbb{S}^{n}$ is defined as
	\begin{equation}
		\operatorname{hvec}(\boldsymbol{W})\triangleq[w_{11}, w_{12}, \ldots, w_{1n}, w_{22}, \ldots, w_{2n},  \ldots,  w_{nn}]^\top,
	\end{equation}
	where $\operatorname{hvec}(\boldsymbol{W})\in\mathbb{R}^{\frac{n^{2}+n}{2}}$.
\end{definition}
	\begin{definition}\label{defi3}
		The matrix operator of any vector $\boldsymbol{x}=[x_{1}, x_{2}, \ldots, x_{n}]^\top\in\mathbb{R}^{n}$ is defined as
		\begin{equation}
			\operatorname{Mat}(\boldsymbol{x})\triangleq
			\left[
			\begin{array}{ccccccccc}
				x_{1} & x_{2}& \ldots& x_{n}&0&\ldots&0&\ldots\ldots&0\\
				0     & x_{1}&\ldots & 0    & x_{2}&\ldots&x_{n}&\ldots\ldots&0\\
				\vdots     &\vdots     &\ddots&\vdots&\vdots&\ddots&\vdots&\ddots&\vdots\\
				0     &0     &\ldots&x_{1} &0&\ldots&x_{2}&\ldots\ldots&x_{n}\\
			\end{array}
			\right],
		\end{equation}
		where 	$\operatorname{Mat}(\boldsymbol{x})\in\mathbb{R}^{n\times\frac{n^{2}+n}{2}}$.
	\end{definition}
	
	\begin{definition}\label{defi2}
		The quadratic mapping operator of any vector $\boldsymbol{x}=[x_{1}, x_{2},  \ldots, x_{n}]^\top\in\mathbb{R}^{n}$ is defined as
		\begin{equation}
			\operatorname{qvec}(\boldsymbol{x})\triangleq[\dfrac{1}{2}x_{1}^{2}, x_{1}x_{2}, \ldots, x_{1}x_{n}, \dfrac{1}{2}x_{2}^{2}, \ldots, x_{2}x_{n}, \ldots, \dfrac{1}{2}x_{n}^{2}]
			^\top,
		\end{equation}
		where $	\operatorname{qvec}(\boldsymbol{x})\in\mathbb{R}^{\frac{n^{2}+n}{2}}$.
	\end{definition}
	
		\subsection{SVM$_{0/1}$}
	SVM$_{0/1}$ \cite{15921057020221001} is a variant of SVM based on the 0-1 loss function. Specifically, for a given training set
	\begin{equation}
		T=\{(\boldsymbol{x}_{i}, y_{i})_{i=1,\dots, N}\mid \boldsymbol{x}_{i}\in\mathbb{R}^{n}, y_{i}\in\{-1, 1\}\}\label{equ:dataset},
	\end{equation}
	where $N$ is the number of samples, $\boldsymbol{x}_{i}\in\mathbb{R}^{n}$ indicates the $i$-th training sample in the $n$-dimensional space, $y_{i}$ is the label of sample $\boldsymbol{x}_{i}$. The task of SVM$_{0/1}$ is to find the following hyperplane
		\begin{equation}
		f(\boldsymbol{x})\triangleq\boldsymbol{w}^\top\boldsymbol{x}+b=0,
	\end{equation}
	where $\boldsymbol{w}\in\mathbb{R}^{n}, b\in\mathbb{R}$.
	Based on the idea of maximum margin and by introducing the 0-1 loss function, the optimization problem is constructed as follows
	\begin{equation}
		\begin{aligned}
			\min_{\boldsymbol{w}, b}\quad&\dfrac{1}{2}\|\boldsymbol{w}\|_{2}^{2}+C\sum_{i=1}^{N}\ell_{0/1}(1-y_{i}f(\boldsymbol{x}_{i})),
			\label{equ01svm}\\
		\end{aligned}
	\end{equation}
	where $\ell_{0/1}(\cdot)$ is the 0-1 loss function, its specific form is
	\begin{equation}
		\ell_{0/1}(1-y_{i} f(\boldsymbol{x}_{i}))=
		\begin{cases}
			1, &\text{$(1-y_{i} f(\boldsymbol{x}_{i}))>0$},\\
			0, &\text{$(1-y_{i} f(\boldsymbol{x}_{i}))\leq0$}, \label{01}
		\end{cases}
	\end{equation}
	and $C>0$ is a penalty parameter.
	
	After acquiring the solution $(\boldsymbol{w}^{*}, b^{*})$ to the optimization problem (\ref{equ01svm}), for a new sample $\boldsymbol{x}\in\mathbb{R}^{n}$, its label is either assigned to $+$ 1 or $-$ 1, by the following decision rule
	\begin{equation}
		g(\boldsymbol{x})=\operatorname{sign}(\boldsymbol{w}^{*\top}\boldsymbol{x}+b^{*}).
	\end{equation}
	
	The optimization problem (\ref{equ01svm}) can be rewritten as an equivalent problem with equality constraint, a fast alternating direction method of multipliers is introduced to solve it. Moreover, the computational complexity of SVM$_{0/1}$ is low.
	
	\subsection{SQSSVM}
	For the training set $T$ (\ref{equ:dataset}), the purpose of SQSSVM \cite{MR3580868} is to solve the following optimization problem,
	\begin{equation}
		\begin{aligned}
			\min_{\boldsymbol{W}, \boldsymbol{b}, c, \boldsymbol{\xi}}\quad
			&\dfrac{1}{2}\sum_{i=1}^{N}\|\boldsymbol{W}\boldsymbol{x}_{i}+\boldsymbol{b}\|_{2}^{2}+C\sum_{i=1}^{N}\xi_{i},\\
			s.t.\quad
			&y_{i}(\dfrac{1}{2}\boldsymbol{x}_{i}^\top\boldsymbol{W}\boldsymbol{x}_{i}+\boldsymbol{b}^\top\boldsymbol{x}_{i}+c)
			\geq1-\xi_{i}, \\
			&\xi_{i}\geq0, \quad i=1, \ldots\ldots, N, 
		\end{aligned}
		\label{sqssvm}
	\end{equation}
	where $C>0$ is a penalty parameter, $\boldsymbol{\xi}=[\xi_{1}, \xi_{2},\ldots, \xi_{N}]^{\top}$ is a non-negative slack vector.
	
	The optimization problem (\ref{sqssvm}) is a convex quadratic programming problem. Both the primal problem and dual problem of SQSSVM can be solved, because they do not involve kernel function.
	
	\section{QSSVM$_{0/1}$}\label{3}
	In this section, a binary classification method is presented, namely QSSVM$_{0/1}$. Then, a fast algorithm and some theoretical characters are proposed for it.
	\subsection{The optimization model of QSSVM$_{0/1}$}\label{3.1jie}
	For the given training set $T$ (\ref{equ:dataset}), we try to find a quadratic separating hyper-surface
	\begin{equation}
		f(\boldsymbol{x})\triangleq\dfrac{1}{2}\boldsymbol{x}^\top \boldsymbol{W}\boldsymbol{x}+\boldsymbol{b}^\top\boldsymbol{x}+c=0,\label{qsurface}
	\end{equation}
	where $\boldsymbol{W}\in \mathbb{S}^{n}, \boldsymbol{b}\in \mathbb{R}^{n}, c\in \mathbb{R}$. Based on the 0-1 loss function and kernel-free trick, the optimization problem is constructed as follows
	\begin{equation}
		\min_{\boldsymbol{W}, \boldsymbol{b}, c} \quad\dfrac{1}{2} \sum_{i=1}^{N}\|\boldsymbol{W} \boldsymbol{x}_{i}+\boldsymbol{b}\|_{2}^{2}+C \sum_{i=1}^{N} \ell_{0/1}(1-y_{i} f(\boldsymbol{x}_{i})),\label{equ01qssvm}
	\end{equation}
	where  $	f(\boldsymbol{x}_{i})\triangleq\dfrac{1}{2}\boldsymbol{x}_{i}^\top \boldsymbol{W}\boldsymbol{x}_{i}+\boldsymbol{b}^\top\boldsymbol{x}_{i}+c=0$, $i=1, \dots, N$, and $C>0$ is a penalty parameter. The 0-1 loss function $\ell_{0/1}(\cdot)$ is defined in the equation~(\ref{01}).
		\begin{figure}[!htbp]
		\centering
		\includegraphics[scale=0.6]{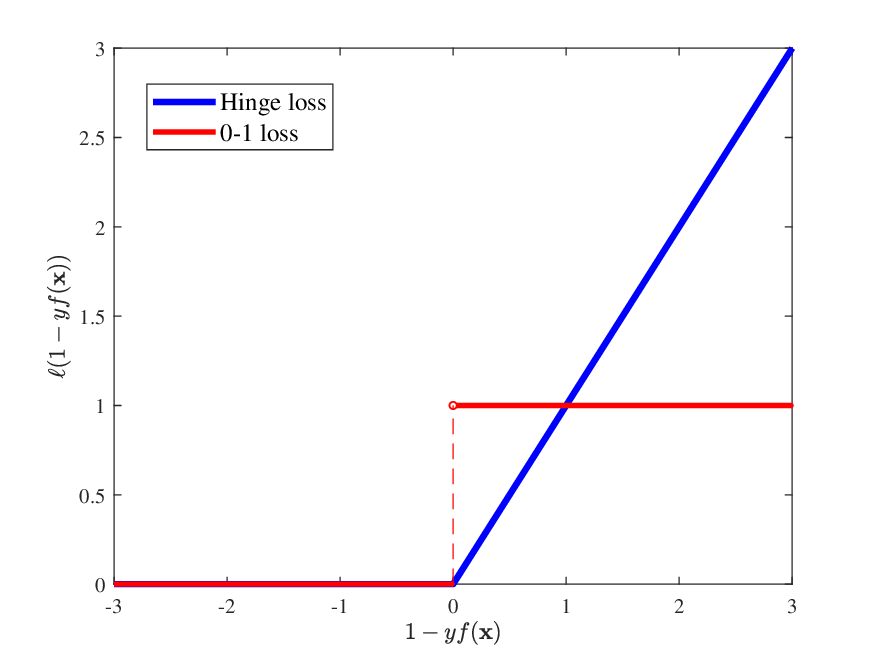}
		\centering
		\caption{The 0-1 loss function and hinge loss function}\label{figure_loss}
	\end{figure}
	 \Cref{figure_loss} shows the geometric interpretability of 0-1 loss function and hinge loss function. They have loss values of 0, when the samples are divided exactly correctly, but behave differently when samples are divided incorrectly or fall within the margin. To be more specific, the loss values of these samples are always 1 for the 0-1 loss function, but increase with the increase of $1-yf(\boldsymbol{x})$ for hinge loss function. Thus, the 0-1 loss function may be more robust compared to the hinge loss function.
	
	In the optimization problem (\ref{equ01qssvm}), the first term of the objective function aims to maximize the sum of relative geometrical margin between all samples and the quadratic hyper-surface, the second term of the objective function is the 0-1 loss function. Let $u_{i}=1-y_{i} f(\boldsymbol{x}_{i})$, then the optimization problem can be written in the form with the equality constraints,
	\begin{equation}
		\begin{aligned}
			\min_{\boldsymbol{W}, \boldsymbol{b}, c, \boldsymbol{u}} \quad&\dfrac{1}{2} \sum_{i=1}^{N}\|\boldsymbol{W} \boldsymbol{x}_{i}+\boldsymbol{b}\|_{2}^
			{2}+C \sum_{i=1}^{N} \ell_{0/1}(u_{i} ),\\
			s.t.\quad& u_{i}+y_{i}f(\boldsymbol{x}_{i})=1, \quad i=1,\ldots, N,\label{yuan1}
		\end{aligned}
	\end{equation}
	where $\boldsymbol{u}=[u_{1}, u_{2}, \ldots, u_{N}]^\top\in\mathbb{R}^{N}$.
	
		Using the symmetry of matrix $\boldsymbol{W}$, let $\widetilde{\boldsymbol{w}}=\operatorname{hvec}(\boldsymbol{W})$ defined by the  \Cref{defi1}, then the above optimization problem can be transformed into the following equivalent form,
	\begin{equation}
		\begin{aligned}
			\min_{\widetilde{\boldsymbol{w}}, \boldsymbol{b}, c, \boldsymbol{u}} \quad&\dfrac{1}{2} \sum_{i=1}^{N}\|\boldsymbol{M}_{i}\widetilde{\boldsymbol{w}} +\boldsymbol{b}\|_{2}^
			{2}+C\|\boldsymbol{u}_{+}\|_{0},\\
			s.t.\quad&\boldsymbol{u}+\boldsymbol{A} \widetilde{\boldsymbol{w}}+\boldsymbol{B}\boldsymbol{b}+c\boldsymbol{y}=\boldsymbol{1},\label{yuan}
		\end{aligned}
	\end{equation}
	where $\boldsymbol{M}_{i}=\operatorname{Mat}(\boldsymbol{x}_{i})\in\mathbb{R}^{n\times\frac{n^{2}+n}{2}}$ defined by the \Cref{defi3},
	$\boldsymbol{u}_{+}=[(u_{1})_{+}, (u_{2})_{+}, \ldots, (u_{N})_{+}]^\top\in\mathbb{R}^{N}$, 
	$(u_{i})_{+}=\max\{u_{i}, 0\}$, 
	$\boldsymbol{1}=[1, 1, \ldots, 1]^\top\in\mathbb{R}^{N}$, $\boldsymbol{y}=[y_{1}, y_{2},\ldots, y_{N}]^\top\in\mathbb{R}^{N}$,  $\boldsymbol{A}=[y_{1}\boldsymbol{s}_{1}, y_{2}\boldsymbol{s}_{2}, \ldots, y_{N}\boldsymbol{s}_{N}]^\top
	\in\mathbb{R}^{N\times \frac{n^{2}+n}{2}}$,
	$\boldsymbol{B}=[y_{1}\boldsymbol{x}_{1}, y_{2}\boldsymbol{x}_{2}, \ldots, y_{N}\boldsymbol{x}_{N}]^\top
	\in\mathbb{R}^{N\times n}$, $\boldsymbol{s}_{i}\triangleq\operatorname{qvec}(\boldsymbol{x}_{i})\in\mathbb{R}^{\frac{n^{2}+n}{2}}$ defined by the \Cref{defi2} and $i=1, \dots, N$.
	
	\subsection{Optimization algorithm}
	In this subsection, a fast algorithm is designed for the optimization problem (\ref{yuan}), which is executed under the ADMM algorithm framework. However, the 0-1 loss function is not differentiable, the P-stationary point of the optimization problem is defined by using the proximal operator of the 0-1 loss function. In this way, the non-convex discontinuity of the optimization problem due to the 0-1 loss function is solved.
	
	\subsubsection{\bf First order optimality condition}
	
	The optimality condition plays a crucial role in solving optimization problem. Next, we propose the proximal stationary point for the optimization problem (\ref{yuan}) and analyze the relationship between it and local/global minimizer.
	
		Now, by introducing a Lagrangian multiplier vector $\boldsymbol{\lambda}\in \mathbb R^{N}$, the Lagrange function of the optimization problem (\ref{yuan}) can be expressed as follows,
	\begin{equation}
		\begin{aligned}
			L(\widetilde{\boldsymbol{w}}; \boldsymbol{b}; c; \boldsymbol{u};
			\boldsymbol{\lambda})=&\dfrac{1}{2} \sum_{i=1}^{N}\|\boldsymbol{M}_{i}\widetilde{\boldsymbol{w}} +\boldsymbol{b}\|_{2}^
			{2}+C\|\boldsymbol{u}_{+}\|_{0}\\
			&+\boldsymbol{\lambda}^{\top}(\boldsymbol{u}+\boldsymbol{A} \widetilde{\boldsymbol{w}}+\boldsymbol{B}\boldsymbol{b}+c\boldsymbol{y}-\boldsymbol{1}).
		\end{aligned}
	\end{equation}
	For vector $\widetilde{\boldsymbol{w}}, \boldsymbol{b}$, scalar $c$, and vector $\boldsymbol{\lambda}$, the following equations can be obtained:
	\begin{eqnarray}
		\nabla_{\widetilde{\boldsymbol{w}}}{L}&=&	\sum_{i=1}^{N}\boldsymbol{M}_{i}^\top(\boldsymbol{M}_{i}\widetilde{\boldsymbol{w}} +\boldsymbol{b})+\boldsymbol{A}^\top\boldsymbol{\lambda},\label{tl1}\\
		\nabla_{\boldsymbol{b}}{L}&=&	\sum_{i=1}^{N}(\boldsymbol{M}_{i}\widetilde{\boldsymbol{w}} +\boldsymbol{b})+\boldsymbol{B}^\top\boldsymbol{\lambda}, \label{tl2}\\
		\nabla_{c}{L}&=&\boldsymbol{y}^\top\boldsymbol{\lambda},\label{tl3}\\
		\nabla_{\boldsymbol{\lambda}}{L}&=&	\boldsymbol{u}+\boldsymbol{A} \widetilde{\boldsymbol{w}}+\boldsymbol{B}\boldsymbol{b}+c\boldsymbol{y}-\boldsymbol{1},\label{tl4}
	\end{eqnarray}
	Since the $\|\boldsymbol{u}_{+}\|_{0}:=\ell_{0/1}(\boldsymbol{u})$, the gradient of $\boldsymbol{u}$ in the Lagrangian function cannot be calculated. Then we obtain the following equation
		\begin{equation}
		\operatorname{prox}_{\gamma C\|(\cdot)_{+}\|_{0}}(\boldsymbol{u}-\gamma\boldsymbol{\lambda})
		=\boldsymbol{u},\label{lt5}
	\end{equation} 
	where
	\begin{equation}
		[\operatorname{prox}_{\gamma C\|(\cdot)_{+}\|_{0}}(\boldsymbol{v})]_{i}=
		\begin{cases}
			0, & 0<v_{i}\leq\sqrt{2 \gamma C}, \\
			v_{i}, & v_{i}>\sqrt{2 \gamma C}\text { or }v_{i}\leq0, 
		\end{cases}
		\label{popetor1}
	\end{equation}
	and $\boldsymbol{v}:=\boldsymbol{u}-\gamma\boldsymbol{\lambda}^{*}$.
	Formula (\ref{popetor1}) is the proximal operator of 0-1 loss function \cite{15921057020221001}. According to the above derivation process, the following definition can be obtained:
	\begin{definition}\label{d1}
		For a given $C>0$, we say $(\widetilde{\boldsymbol{w}}^{*}; \boldsymbol{b}^{*}; c^{*}; \boldsymbol{u}^{*})$ is a proximal stationary (P-stationary) point of the optimization problem (\ref{yuan}), if there exists a Lagrangian multiplier vector $\boldsymbol{\lambda}^{*}\in \mathbb R^{N}$, and a constant $\gamma>0$ such that 
		\begin{eqnarray}
			\sum_{i=1}^{N}\boldsymbol{M}_{i}^\top(\boldsymbol{M}_{i}\widetilde{\boldsymbol{w}}^{*} +\boldsymbol{b}^{*})+\boldsymbol{A}^\top\boldsymbol{\lambda}^{*}&=&\boldsymbol{0},\label{l1}\\
			\sum_{i=1}^{N}(\boldsymbol{M}_{i}\widetilde{\boldsymbol{w}}^{*} +\boldsymbol{b}^{*})+\boldsymbol{B}^\top\boldsymbol{\lambda}^{*}&=&\boldsymbol{0}, \label{l2}\\
			\boldsymbol{y}^\top\boldsymbol{\lambda}^{*}&=&0,\label{l3}\\
			\boldsymbol{u}^{*}+\boldsymbol{A} \widetilde{\boldsymbol{w}}^{*}+\boldsymbol{B}\boldsymbol{b}^{*}+c\boldsymbol{y}-\boldsymbol{1}&=&\boldsymbol{0},\label{l4}\\
			\operatorname{prox}_{\gamma C\|(\cdot)_{+}\|_{0}}(\boldsymbol{u}^{*}-\gamma\boldsymbol{\lambda}^{*})
			&=&\boldsymbol{u}^{*},\label{l5}
		\end{eqnarray}
	\end{definition}
	
		Further exploration shows that the P-stationary point is closely related to the local or global minimizer of the optimization problem (\ref{yuan}). Let
	\begin{eqnarray}
		\boldsymbol{D}&:=&[\boldsymbol{A}\, \boldsymbol{B}\, \boldsymbol{y}]\in\mathbb{R}^{N\times (\frac{n^{2}+n}{2}+n+1)}, \label{D}\\
		\boldsymbol{E}_{i}&:=&
		[\boldsymbol{M}_{i}\,
		\boldsymbol{I}_{n}\,  \boldsymbol{0}_{n}
		]\in\mathbb{R}^{n\times(\frac{n^{2}+n}{2}+n+1)},\label{Ei}\\
		\boldsymbol{H}_{i}&:=&\boldsymbol{E}_{i}\boldsymbol{D}^{+}\in\mathbb{R}^{n\times N},\label{Hi}\\
		\boldsymbol{H}&:=&\sum_{i=1}^{N}\boldsymbol{H}_{i}^\top\boldsymbol{H}_{i}\in\mathbb{R}^{N\times N},\label{H}
	\end{eqnarray}
	where the matrix $\boldsymbol{D}^{+}\in\mathbb{R}^{(\frac{n^{2}+n}{2}+n+1)\times N}$ is the generalized inverse of matrix $\boldsymbol{D}$, $i=1, \dots, N$. Then we have the following theorem.
		\begin{theorem} There are the following relationships in optimization problem (\ref{yuan}), \label{theorem}
		
		(1) For a given $C>0$, if $\boldsymbol{D}$ is a full column rank matrix, such that the globally optimal solution is a P-stationary point with $0<\gamma<\lambda_{\max}(\boldsymbol{H})$.
		
		(2) For a given $C>0$, if $(\widetilde{\boldsymbol{w}}^{*}; \boldsymbol{b}^{*}; c^{*}; \boldsymbol{u}^{*})$ is a P-stationary point with $\gamma>0$, which is also a locally optimal solution.
	\end{theorem}
	The above Theorem \ref{theorem} summarizes the optimality conditions of QSSVM$_{0/1}$, and its proof is shown in \ref{zhengming1}.1.

	\subsubsection{\bf Support vectors of QSSVM$_{0/1}$}
	
	Inspired by the significant role of support vectors in the training process, the support vectors of QSSVM$_{0/1}$ are defined as follows.
	\begin{definition}
		Suppose that $\boldsymbol{\lambda}^{*}=(\lambda_{1}^{*}, \lambda_{2}^{*}, \ldots, \lambda_{N}^{*})$ is a Lagrangian multiplier vector of QSSVM$_{0/1}$, the input vector $\boldsymbol{x}_{i}$ of the training sample $(\boldsymbol{x}_{i}, y_{i})$ is called a support vector, if the corresponding component $\lambda_{i}^{*}$ of $\boldsymbol{\lambda}^{*}$ is nonzero, and otherwise it is a nonsupport vector.
	\end{definition}
	The support vectors of QSSVM$_{0/1}$ are derived as follows.
	
		Firstly, let $(\widetilde{\boldsymbol{w}}^{*}; \boldsymbol{b}^{*}; c^{*}; \boldsymbol{u}^{*})$ be a P-stationary point, thus it is a local minimizer of the optimization problem (\ref{yuan}).
	According to the range of values of the variables in the formula~(\ref{popetor1}) , let
	\begin{equation}
		T_{*}:=\{i \in \mathbb{N}_{N}:u_{i}^{*}-\gamma\lambda_{i}^{*} \in(0, \sqrt{2 \gamma C}]\},\label{T*}
	\end{equation}
	and $\overline{{T}}_{*}:=\mathbb{N}_{N}\setminus T_{*}$ be a complement of $T_{*}$.
	So, the formula (\ref{l5}) can be represented as
	\begin{equation}
		\left[\begin{array}{c}
			(\operatorname{prox}_{\gamma C\|(\cdot)_{+}\|_{0}}	(\boldsymbol{u}^{*}-\gamma \boldsymbol{\lambda}^{*}))_{T_{*}} \\
			(\operatorname{prox}_{\gamma C\|(\cdot)_{+}\|_{0}}(\boldsymbol{u}^{*}-\gamma \boldsymbol{\lambda}^{*}))_{\overline{T}_{*}}
		\end{array}\right]
		=\left[\begin{array}{c}\boldsymbol{0}_{T_{*}} \\
			(\boldsymbol{u}^{*}-\gamma \boldsymbol{\lambda}^{*})_{\overline{T}_{*}}
		\end{array}\right]. \label{ut*}
	\end{equation}
	From this, we can obtain
	\begin{equation}
		\left[\begin{array}{c}
			\boldsymbol{u}_{T_{*}}^{*} \\
			\boldsymbol{\lambda}^{*}_{\overline{T}_{*}}
		\end{array}\right]=\boldsymbol{0}.\label{uandlamda}
	\end{equation}
	Due to that $u_i{}^{*}=0$ in the index set $T_{*}$,
	and $\lambda_{i}^{*}=0$ in the index set $\overline{T}_{*}$,
	the value of $\lambda_{i}^{*}$ is
	\begin{equation}
		\lambda_{i}^{*} \begin{cases}\in[-\sqrt{2C/ \gamma}, 0), & i\in T_{*}, \\=0, & i\in\overline{T}_{*}.\end{cases}\label{lamdaq}
	\end{equation}
	
	Secondly, put the formula (\ref{lamdaq}) into (\ref{l1}) and (\ref{l2}), then the solution to the optimization problem (\ref{yuan}) with respect to ($\widetilde{\boldsymbol{w}}^{*}; \boldsymbol{b}^{*}$)  is indicated as
	\begin{equation}
		\begin{aligned}
			[\widetilde{\boldsymbol{w}}^{*}; \boldsymbol{b}^{*}]&=-(\sum_{i=1}^{N}[\boldsymbol{M}_{i}, \boldsymbol{I}_{n}]^\top[\boldsymbol{M}_{i}, \boldsymbol{I}_{n}])^{-1}[\boldsymbol{A}_{T_{*}}, \boldsymbol{B}_{T_{*}}]^\top\boldsymbol{\lambda}_{T_{*}}^{*}\\
			&=-\sum_{i\in T_{*}}\lambda_{i}^{*}y_{i}\boldsymbol{G}^{-1}[\boldsymbol{s}_{i}^{\top}, \boldsymbol{x}_{i}^{\top}]^{\top},
		\end{aligned}\label{zdelagelangri}
	\end{equation}
	where $\boldsymbol{G}=\sum_{i=1}^{N}[\boldsymbol{M}_{i}, \boldsymbol{I}_{n}]^\top[\boldsymbol{M}_{i}, \boldsymbol{I}_{n}]$.
	Taking $\boldsymbol{u}_{T_{*}}^{*}=\boldsymbol{0}$ into the formula (\ref{l4}) to obtain
	\begin{equation}
		(\boldsymbol{A} \widetilde{\boldsymbol{w}}^{*}+\boldsymbol{B}\boldsymbol{b}^{*}+c^{*}\boldsymbol{y})_{T_{*}}=\boldsymbol{1},
	\end{equation}
	namely,
	\begin{equation}
		\boldsymbol{s}_{i}^\top \widetilde{\boldsymbol{w}}^{*}+\boldsymbol{b}^{*\top}\boldsymbol{x}_{i}+c^{*}
		=\pm1, \quad(i\in T_{*}). \label{zhichengqumian}
	\end{equation}
	The above equation is equivalent to \begin{equation}\dfrac{1}{2}\boldsymbol{x}_{i}^\top \boldsymbol{W}^{*}\boldsymbol{x}_{i}+\boldsymbol{b}^{*\top}\boldsymbol{x}_{i}+c^{*}=\pm1, \quad(i\in T_{*}).\label{hyper-surface}
	\end{equation}			
	In fact, the Lagrangian multiplier $-\boldsymbol{\lambda}^{*}$ is a solution to the dual problem of the optimization problem (\ref{yuan}).
	Moreover, $\{\boldsymbol{x}_{i}, i\in T_{*}\}$ are standard support vectors (SVs), which are selected by the formula (\ref{ut*}). More importantly, the formula (\ref{hyper-surface}) shows that all SVs fall on the support hyper-surfaces. 
	This indicates that the SVs arrange orderly and may also lead to sparsity in the model.
	
\subsubsection{\bf ADMM algorithm of QSSVM$_{0/1}$}\label{3.22}
	In order to solve the optimization problem (\ref{yuan}), a new iterative algorithm based on ADMM algorithm framework is proposed.
	
		The augmented Lagrangian function of the optimization problem (\ref{yuan}) is given by
	\begin{equation}
		\begin{aligned}
			L_{\sigma}(\widetilde{\boldsymbol{w}}; \boldsymbol{b}; c; \boldsymbol{u};
			\boldsymbol{\lambda})=&\dfrac{1}{2} \sum_{i=1}^{N}\|\boldsymbol{M}_{i}\widetilde{\boldsymbol{w}} +\boldsymbol{b}\|_{2}^
			{2}+C\|\boldsymbol{u}_{+}\|_{0}\\
			&+\boldsymbol{\lambda}^{\top}(\boldsymbol{u}-\boldsymbol{1}+\boldsymbol{A} \widetilde{\boldsymbol{w}}+\boldsymbol{B}\boldsymbol{b}+c\boldsymbol{y})\\
			&+\dfrac{\sigma}{2}\|\boldsymbol{u}-\boldsymbol{1}+\boldsymbol{A} \widetilde{\boldsymbol{w}}+\boldsymbol{B}\boldsymbol{b}+c\boldsymbol{y}\|^{2},
		\end{aligned}
	\end{equation}
	where $\boldsymbol{\lambda}$ is a Lagrangian multiplier vector, $\sigma>0$ is a penalty parameter. Given the $k$-th iteration $(\widetilde{\boldsymbol{w}}^{k}; \boldsymbol{b}^{k}; c^{k}; \boldsymbol{u}^{k}; \boldsymbol{\lambda}^{k})$, the solution framework is as follows
	\begin{equation}
		\begin{aligned}
			&\boldsymbol{u}^{k+1}
			=\arg\min_{\boldsymbol{u}}L_{\sigma}(\widetilde{\boldsymbol{w}}^{k}; \boldsymbol{b}^{k}; c^{k}; \boldsymbol{u}; \boldsymbol{\lambda}^{k}),\\
			&\begin{aligned}
				[\widetilde{\boldsymbol{w}}^{k+1}; \boldsymbol{b}^{k+1}]
				=&\arg\min_{[\widetilde{\boldsymbol{w}}; \boldsymbol{b}]}L_{\sigma}([\widetilde{\boldsymbol{w}}; \boldsymbol{b}]; c^{k}; \boldsymbol{u}^{k+1}; \boldsymbol{\lambda}^{k})\\
				&	+
				\dfrac{\sigma}{2}\|[\widetilde{\boldsymbol{w}}; \boldsymbol{b}]-[\widetilde{\boldsymbol{w}}^{k}; \boldsymbol{b}^{k}]\|_{D_{k}}^{2},
			\end{aligned}\\
			&c^{k+1}
			=\arg\min_{c}L_{\sigma}(\widetilde{\boldsymbol{w}}^{k+1}; \boldsymbol{b}^{k+1}; c; \boldsymbol{u}^{k+1}; \boldsymbol{\lambda}^{k}),\\
			&\boldsymbol{\lambda}^{k+1}
			=\boldsymbol{\lambda}^{k}+\eta\sigma(\boldsymbol{u}^{k+1}-\boldsymbol{1}
			+\boldsymbol{A}\widetilde{\boldsymbol{w}}^{k+1}	+\boldsymbol{B}\boldsymbol{b}^{k+1}+c^{k+1}\boldsymbol{y}), \\
		\end{aligned}\label{update}
	\end{equation}
	where $\eta>0$ is known as the dual step-size. The proximal 
	term with respect to  $\widetilde{\boldsymbol{w}}$ and $\boldsymbol{b}$ is 
	\begin{equation}
		\begin{aligned}
			&\|[\widetilde{\boldsymbol{w}}; \boldsymbol{b}]-[\widetilde{\boldsymbol{w}}^{k}; \boldsymbol{b}^{k}] \|_{D_{k}}^{2}\\
			=&\langle[\widetilde{\boldsymbol{w}}; \boldsymbol{b}]-[\widetilde{\boldsymbol{w}}^{k}; \boldsymbol{b}^{k}], D_{k}([\widetilde{\boldsymbol{w}}; \boldsymbol{b}] -[\widetilde{\boldsymbol{w}}^{k}; \boldsymbol{b}^{k}])\rangle,\label{Dk}
		\end{aligned}
	\end{equation}
	where $D_{k}$
	is a symmetric matrix. 
	Each subproblem in the formula (\ref{update}) is updated as follows.\\
	$\boldsymbol{step1}$. Updating $\boldsymbol{u}^{k+1}$:
	\begin{equation}
		\begin{aligned}
&\boldsymbol{u}^{k+1}\\
			&\begin{aligned}=&\arg\min_{\boldsymbol{u} }C\|\boldsymbol{u}_{+}\|_{0}+\boldsymbol{\lambda}^{k\top}\boldsymbol{u}\\
				&+\dfrac{\sigma}{2}\|\boldsymbol{u}-\boldsymbol{1}+\boldsymbol{A} \widetilde{\boldsymbol{w}}^{k}+\boldsymbol{B}\boldsymbol{b}^{k}+c^{k}\boldsymbol{y}\|^{2}\end{aligned}\\
			&=\arg\min_{\boldsymbol{u} }C\|\boldsymbol{u}_{+}\|_{0}+\dfrac{\sigma}{2}\|\boldsymbol{u}-\boldsymbol{v}^{k}\|^{2}\\
			&=\operatorname{prox}_{\frac{1}{\sigma} C\|(\cdot)_{+}\|_{0}}(\boldsymbol{v}^{k}).\label{uk+1}
		\end{aligned}
	\end{equation}
	
		where the third equation is derived from the Definition 3.2 \cite{15921057020221001} with $\sigma=1/\gamma$, and 
	$\boldsymbol{v}^{k}:=\boldsymbol{1}-\boldsymbol{A} \widetilde{\boldsymbol{w}}^{k}-\boldsymbol{B}\boldsymbol{b}^{k}-c^{k}\boldsymbol{y}-{\boldsymbol{\lambda}^{k}}/{\sigma}$. According to $\boldsymbol{v}^{k}$, define a working set
	\begin{equation}
		T_{k}=\{i\in\mathbb{N}_{N}: v_{i}^{k}\in(0, \sqrt{2C/\sigma}]\}, \label{updateTk+1}
	\end{equation}
	and its complement $\overline{T}_{k}:=\mathbb{N}_{N}\setminus T_{k}$, we have
	\begin{equation}
		\boldsymbol{u}^{k+1}_{T_{k}}=\boldsymbol{0},\,\,
		\boldsymbol{u}^{k+1}_{\overline{T}_{k}}=\boldsymbol{v}^{k}_{\overline{T}_{k}}.\label{updateuk+1}
	\end{equation}
	$\boldsymbol{step2}$. Updating $[\widetilde{\boldsymbol{w}}^{k+1}; \boldsymbol{b}^{k+1}]$: the subproblem with respect to $[\widetilde{\boldsymbol{w}}; \boldsymbol{b}]$ in the formula (\ref{update}) is
	\begin{equation}
		\begin{aligned}	
			[\widetilde{\boldsymbol{w}}^{k+1}; \boldsymbol{b}^{k+1}]=&\arg\min_{[\widetilde{\boldsymbol{w}}; \boldsymbol{b}]}
			\dfrac{1}{2} \sum_{i=1}^{N}\|[\boldsymbol{M}_{i}\:\boldsymbol{I}_{n}][\widetilde{\boldsymbol{w}}; \boldsymbol{b}]\|_{2}^
			{2}\\
			&+\boldsymbol{\lambda}^{k\top}[\boldsymbol{A}\:\boldsymbol{B}][\widetilde{\boldsymbol{w}}; \boldsymbol{b}]+\dfrac{\sigma}{2}\|[\widetilde{\boldsymbol{w}}; \boldsymbol{b}]-[\widetilde{\boldsymbol{w}}^{k}; \boldsymbol{b}^{k}]\|_{D_{k}}^{2}\\
			&+\dfrac{\sigma}{2}\|\boldsymbol{u}^{k+1}-\boldsymbol{1}+[\boldsymbol{A}\:\boldsymbol{B}] [\widetilde{\boldsymbol{w}}; \boldsymbol{b}]
			+c^{k}\boldsymbol{y}\|^{2},
		\end{aligned}
	\end{equation}
	where $D_{k}=-[\boldsymbol{A}_{\overline{T}_{k}}\:\boldsymbol{B}_{\overline{T}_{k}}]^{\top} [\boldsymbol{A}_{\overline{T}_{k}}\:\boldsymbol{B}_{\overline{T}_{k}}]$.
	Setting the gradient of $L_{\sigma}([\widetilde{\boldsymbol{w}}; \boldsymbol{b}])$ with respect to $[\widetilde{\boldsymbol{w}}; \boldsymbol{b}]$ to be zero, we have
	\begin{equation}
		\begin{aligned}
			&\sum_{i=1}^{N}[\boldsymbol{M}_{i}\:\boldsymbol{I}_{n}]^{\top}[\boldsymbol{M}_{i}\:\boldsymbol{I}_{n}][\widetilde{\boldsymbol{w}}; \boldsymbol{b}]+[\boldsymbol{A}\:\boldsymbol{B}]^{\top} \boldsymbol{\lambda}^{k}\\
			& +\sigma [\boldsymbol{A}\:\boldsymbol{B}] ^{\top}(\boldsymbol{u}^{k+1}-\boldsymbol{1}+[\boldsymbol{A}\:\boldsymbol{B}] [\widetilde{\boldsymbol{w}}; \boldsymbol{b}]
			+c^{k}\boldsymbol{y})\\
			&-\sigma [\boldsymbol{A}_{\overline{T}_{k}}\:\boldsymbol{B}_{\overline{T}_{k}}]^{\top} [\boldsymbol{A}_{\overline{T}_{k}}\:\boldsymbol{B}_{\overline{T}_{k}}]([\widetilde{\boldsymbol{w}}; \boldsymbol{b}]-[\widetilde{\boldsymbol{w}}^{k}; \boldsymbol{b}^{k}])=\boldsymbol{0}.
		\end{aligned}
	\end{equation}
		Since \begin{equation}\begin{aligned}
			&[\boldsymbol{A}_{T_{k}}\:\boldsymbol{B}_{T_{k}}]^{\top}[\boldsymbol{A}_{T_{k}}\:\boldsymbol{B}_{T_{k}}]\\
			=&[\boldsymbol{A}\:\boldsymbol{B}]^{\top} [\boldsymbol{A}\:\boldsymbol{B}]-[\boldsymbol{A}_{\overline{T}_{k}}\:\boldsymbol{B}_{\overline{T}_{k}}]^{\top} [\boldsymbol{A}_{\overline{T}_{k}}\:\boldsymbol{B}_{\overline{T}_{k}}],
		\end{aligned}
	\end{equation}
	it is equivalent to find a solution to the linear system of equations
	\begin{equation}
		(\boldsymbol{G}+ \sigma [\boldsymbol{A}_{T_{k}}\:\boldsymbol{B}_{T_{k}}]^{\top}[\boldsymbol{A}_{T_{k}}\:\boldsymbol{B}_{T_{k}}])[\widetilde{\boldsymbol{w}}; \boldsymbol{b}]
		=\sigma [\boldsymbol{A}_{T_{k}}\:\boldsymbol{B}_{T_{k}}]^{\top}\boldsymbol{d}^{k}
		\label{updatew},
	\end{equation}
	where $\boldsymbol{G}=\sum_{i=1}^{N}[\boldsymbol{M}_{i}, \boldsymbol{I}_{n}]^\top[\boldsymbol{M}_{i}, \boldsymbol{I}_{n}]$,  $\boldsymbol{d}^{k}:=-(\boldsymbol{u}^{k+1}+c^{k}\boldsymbol{y}-\boldsymbol{1}+\boldsymbol{\lambda}^{k}/\sigma)$.
	
		Then $[\widetilde{\boldsymbol{w}}^{k+1}; \boldsymbol{b}^{k+1}]$ can be updated in the following ways,
	\begin{itemize}
		\item if $({n^{2}+3n})/{2}\leq|T_{k}|$, the linear system of equations (\ref{updatew}) can be directly solved by 
		\begin{equation}
			\begin{aligned}
				&	[\widetilde{\boldsymbol{w}}^{k+1}; \boldsymbol{b}]\\
				=&(\boldsymbol{G}+\sigma[\boldsymbol{A}_{T_{k}}\:\boldsymbol{B}_{T_{k}}]^{\top}[\boldsymbol{A}_{T_{k}}\:\boldsymbol{B}_{T_{k}}])^{-1}(\sigma [\boldsymbol{A}_{T_{k}}\:\boldsymbol{B}_{T_{k}}]^{\top} \boldsymbol{d}_{T_{k}}^{k}).\label{updatezk+1}
			\end{aligned}
		\end{equation}
		\item if $({n^{2}+3n})/{2}>|T_{k}|$, the conjugate gradient (CG) method is used to solve the linear system of equations (\ref{updatew}) for efficiency.
		\end{itemize}
	$\boldsymbol{step3}$. Updating $c^{k+1}$: the subproblem with respect to $c$ in the formula (\ref{update}) is
	\begin{equation}
		c^{k+1}=\arg\min_{c}
		\: \boldsymbol{\lambda}^{k\top}c\boldsymbol{y}
		+\dfrac{\sigma}{2}\|\boldsymbol{u}^{k+1}-\boldsymbol{1}+\boldsymbol{A}\widetilde{\boldsymbol{w}}^{k+1}+\boldsymbol{B}\boldsymbol{b}^{k+1}+c\boldsymbol{y}\|^{2}.
	\end{equation}
	Setting the derivation of $L_{\sigma}(c)$ with respect to $c$ to be zero, we have
	\begin{equation}
		\boldsymbol{y}^{\top}\boldsymbol{\lambda}^{k}+\sigma\boldsymbol{y}^{\top}(\boldsymbol{u}^{k+1}-\boldsymbol{1}+
		\boldsymbol{A}\widetilde{\boldsymbol{w}}^{k+1}+\boldsymbol{B}\boldsymbol{b}^{k+1}+c\boldsymbol{y})=0.
	\end{equation}
	Then $c^{k+1}$ can be updated by 
		\begin{equation}
		c^{k+1}
		=-\boldsymbol{y}^{\top}(\boldsymbol{A}\widetilde{\boldsymbol{w}}^{k+1}+\boldsymbol{B} \boldsymbol{b}^{k+1}-\boldsymbol{1}+\boldsymbol{u}^{k+1}+\boldsymbol{\lambda}^{k}/\sigma)/N.\label{updateck+1}
	\end{equation}
	$\boldsymbol{step4}$. Updating $\boldsymbol{\lambda}^{k+1}:$
	\begin{equation}
		\begin{aligned}
			&\boldsymbol{\lambda}_{T_{k}}^{k+1}=\boldsymbol{\lambda}_{T_{k}}^{k}+\eta\sigma(\boldsymbol{u}^{k+1}-\boldsymbol{1}+\boldsymbol{A}\widetilde{\boldsymbol{w}}^{k+1}+\boldsymbol{B}\boldsymbol{b}^{k+1}
			+c^{k+1}\boldsymbol{y}),\\ &\boldsymbol{\lambda}_{\overline{T}_{k}}^{k+1}=\boldsymbol{0}.\label{lamda}
		\end{aligned}
	\end{equation}
	
		Steps 1-4 lead to \Cref{algorithm1} for solving the optimization problem (\ref{yuan}).
	
		\begin{algorithm}[H]
		\caption{\quad ADMM Algorithm } \label{algorithm1}
		\hspace*{0.02in} {\bf Input:} 
		Training set $T$ (\ref{equ:dataset}), parameters $C$, $\sigma$ and $\eta$,  $\boldsymbol{M}_{i}, \: i=1,\ldots,N $,
		the maximum number of iteration steps $K$.\\
		\hspace*{0.02in} {\bf Output:} 
		The final solution $(\widetilde{\boldsymbol{w}}^{k}, \boldsymbol{b}^{k}, c^{k}, \boldsymbol{u}^{k}, \boldsymbol{\lambda}^{k})$ to the optimization problem (\ref{yuan}) .
		\begin{algorithmic}
		\State Initialize($\widetilde{\boldsymbol{w}}^{0}; \boldsymbol{b}^{0}; c^{0}; \boldsymbol{u}^{0}; \boldsymbol{\lambda}^{0}$), set $k=0$.
		\While{The halting condition does not hold and $k\leq K$ } 
		\State Update $T_{k}$ as in (\ref{updateTk+1}).
		\State Update $\boldsymbol{u}^{k+1}$ by formula (\ref{updateuk+1}).
		\If{$({n^{2}+3n})/{2}\leq|T_{k}|$}
		\State Update $[\widetilde{\boldsymbol{w}}^{k+1}; \boldsymbol{b}^{k+1}]$ by formula (\ref{updatezk+1}).
		\Else
		\State Update $[\widetilde{\boldsymbol{w}}^{k+1}; \boldsymbol{b}^{k+1}]$ by the CG method.
		\EndIf
		\State Update $c^{k+1}$ by formula (\ref{updateck+1}).
		\State Update $\boldsymbol{\lambda}^{k+1}$ by formula (\ref{lamda}).
		\State set $k=k+1$.
		\EndWhile
	\end{algorithmic}
	\end{algorithm}
	
	\subsection{The interpretability}
	After obtaining the solution to the optimization problem (\ref{yuan}) with respect to $(\widetilde{\boldsymbol{w}}, \boldsymbol{b}, c)$, the half-vectorization operator defined by Definition \ref{defi1} can be used to invert the symmetric matrix $\boldsymbol{W}^{*}$. 
	Consequently, for a unlabeled sample $\boldsymbol{x}\in\mathbb{R}^{n}$, its label is determined by the following decision rule
	\begin{equation}
		g(\boldsymbol{x})=\operatorname{sign}(\dfrac{1}{2}\boldsymbol{x}^\top\boldsymbol{W}^{*}\boldsymbol{x}+\boldsymbol{b}^{*\top}\boldsymbol{x}+c^{*}).\label{decision}
	\end{equation}
	
	Suppose that $(\boldsymbol{W}^{*}, \boldsymbol{b}^{*}, c^{*})$ in formula (\ref{decision}) is the optimal solution to optimization problem (\ref{equ01qssvm}), then the quadratic hyper-surface (\ref{qsurface}) can be written as
	\begin{equation}
		\begin{aligned}
			f(\boldsymbol{x})&=\dfrac{1}{2}\boldsymbol{x}^\top\boldsymbol{\boldsymbol{W}^{*}}\boldsymbol{x}+\boldsymbol{b^{*}}^\top\boldsymbol{x}+c^{*}\\
			&=\dfrac{1}{2}\sum_{i=1}^{n}\sum_{j=1}^{n}w_{ij}^{*}x_{i}x_{j}+\sum_{i=1}^{n}b_{i}^{*}x_{i}+c^{*}\\
			&=0,\quad i=1, \cdots, m,
		\end{aligned}
	\end{equation}
	where $x_{i}$ is the $i$-th component of the vector $\boldsymbol{x}$, $w_{ij}^{*}$ is the $i$-th row and the $j$-th column
	component of the matrix $\boldsymbol{W}^{*}\in \mathbb{S}^{n}$,
	and $b_{i}^{*}$ is the $i$-th component of the vector $\boldsymbol{b^{*}}\in\mathbb{R}^{n}$.
	Each component of vector $\boldsymbol{x}$ contributes to the quadratic polynomial function. Specifically speaking, when $i=j$, a larger value of $|w_{ii}^{*}|+|b_{i}^{*}|$ indicates a larger contribution from the $i$-th component of the vector $\boldsymbol{x}$.
	When $i\neq j$, the $i$-th and $j$-th components of the vector $\boldsymbol{x}$ act synergistically on the quadratic term.
	Notably, when $\boldsymbol{W}^{*}$ is a zero matrix, the quadratic hyper-surface degenerates into a linear hyperplane. Therefore, our QSSVM$_{0/1}$ is more interpretable than methods with kernel functions.
	
		\subsection{Convergence and complexity analysis of QSSVM$_{0/1}$}
	Since we propose a kernel-free nonlinear method, it avoids choosing the kernel function and tuning the corresponding parameters. But when the sample dimensions are large, it may take more time than SVM$_{0/1}$.  Nevertheless, it still obtains a convergence theorem of Algorithm \ref{algorithm1}. Moreover, the complexity of Algorithm \ref{algorithm1} is analyzed in this subsection.
	\begin{theorem}\label{theorem2}
		Suppose that $(\widetilde{\boldsymbol{w}}^{*}; \boldsymbol{b}^{*}; c^{*}; \boldsymbol{u}^{*}; \boldsymbol{\lambda}^{*})$ is the limit point of the sequence $\{(\widetilde{\boldsymbol{w}}^{k}; \boldsymbol{b}^{k}; c^{k}; \boldsymbol{u}^{k}; \boldsymbol{\lambda}^{k})\}$ generated by Algorithm \ref{algorithm1}, then $(\widetilde{\boldsymbol{w}}^{*}; \boldsymbol{b}^{*}; c^{*}; \boldsymbol{u}^{*})$ is a P-stationary point with $\gamma=1/\sigma$ and also a locally optimal solution to the optimization problem (\ref{yuan}).
	\end{theorem}
	The proof of Theorem \ref{theorem2} is given in \ref{zhengming2}.2. 
	The above Theorem \ref{theorem2} states that if the sequence generated by Algorithm \ref{algorithm1} has a limit point, then it is a P-stationary point and a locally optimal solution to optimization problem (\ref{yuan}). 
	Since the optimization problems (\ref{yuan1}) and (\ref{yuan}) are equivalent, the following theorem are obtained,
	\begin{theorem}
		Suppose that $(\widetilde{\boldsymbol{w}}^{*}; \boldsymbol{b}^{*}; c^{*}; \boldsymbol{u}^{*})$ is a P-stationary point of the optimization problem (\ref{yuan}), then $(\boldsymbol{W}^{*}; \boldsymbol{b}^{*}; c^{*}; \boldsymbol{u}^{*})$ is a P-stationary point of the optimization problem (\ref{yuan1}).		
	\end{theorem}
	
	The computational complexity of each iteration of Algorithm \ref{algorithm1} is analyzed as follows.
	\begin{itemize}
		\item Updating $T_{k}$ by formula (\ref{updateTk+1}) needs the complexity $\mathcal{O}(N)$.
		\item Updating $\boldsymbol{u}^{k+1}$ by formula (\ref{updateuk+1}), the complexity of its main item $(\boldsymbol{A}\widetilde{\boldsymbol{w}}^{k}+\boldsymbol{B}\boldsymbol{b}^{k})$ is $\mathcal{O}(\dfrac{Nn^{2}+3Nn}{2})$.
		\item Updating $[\widetilde{\boldsymbol{w}}^{k+1}; \boldsymbol{b}^{k+1}]$ by formula (\ref{updatezk+1})	if {$({n^{2}+3n})/{2}\leq|T_{k}|$} and CG method otherwise.  For the former, the main computations from  $$[\boldsymbol{A}_{T_{k}}\:\boldsymbol{B}_{T_{k}}]^{\top}[\boldsymbol{A}_{T_{k}}\:\boldsymbol{B}_{T_{k}}], 
		(\boldsymbol{G}+\sigma[\boldsymbol{A}_{T_{k}}\:\boldsymbol{B}_{T_{k}}]^{\top}[\boldsymbol{A}_{T_{k}}\:\boldsymbol{B}_{T_{k}}])^{-1}$$
		with the computational complexities 
		$\mathcal{O}((\dfrac{n^{2}+3n}{2})^{2}|T_{k}|)$ and 	$\mathcal{O}((\dfrac{n^{2}+3n}{2})^{\iota})$ with ${\iota}\in (2,3)$, respectively. For the latter, the computational complexity of calculating  $[\boldsymbol{A}_{T_{k}}\:\boldsymbol{B}_{T_{k}}]^{\top}[\boldsymbol{A}_{T_{k}}\:\boldsymbol{B}_{T_{k}}]$ remains the same as previous, and the computational complexity of the linear system of equations (\ref{updatew}) calculated through the CG method is $\mathcal{O}((\dfrac{n^{2}+3n}{2})^{2}q)$, where $q$ is the number of distinct eigenvalues of $(\boldsymbol{G}+\sigma[\boldsymbol{A}_{T_{k}}\:\boldsymbol{B}_{T_{k}}]^{\top}[\boldsymbol{A}_{T_{k}}\:\boldsymbol{B}_{T_{k}}])$.
		\item As for $c^{k+1}$, which is updated by the formula (\ref{updateck+1}), the complexity of its main item  $(\boldsymbol{A}\widetilde{\boldsymbol{w}}^{k+1}+\boldsymbol{B}\boldsymbol{b}^{k+1})$ is $\mathcal{O}(\dfrac{Nn^{2}+3Nn}{2})$.
		\item Similar to $c^{k+1}$, updating $\boldsymbol{\lambda}^{k+1}$ by the formula (\ref{lamda}), its complexity is $\mathcal{O}(\dfrac{Nn^{2}+3Nn}{2})$.
	\end{itemize}
	Therefore, the whole computational complexity  of the Algorithm \ref{algorithm1} is $\mathcal{O}(\dfrac{Nn^{2}+3Nn}{2}+(\dfrac{n^{2}+3n}{2})^{2}\max\{| T_{k}|, q \}
	)$.
	
	\section{Numerical experiments}\label{4}
	In this section, some experiments are used to verify the performance of our QSSVM$_{0/1}$. It is compared with other state-of-the-art classifiers, including SVM$_{0/1}$ and KSVM$_{0/1}$ using 0-1 loss function, and SVM classifiers with various loss functions, such as the hinge soft-margin loss (SVM), $\nu$SVM \cite{2000New, 0Support} with the hinge soft-margin loss ($\nu$SVM), pinball soft-margin loss (PSVM), square soft-margin loss (LSSVM), ramp soft-margin loss (RSVM) and non-convex robust and smooth soft-margin loss (RSHSVM). Three different kernel functions are used, including linear, RBF and polynomial kernel functions. In addition, the QSSVM$_{0/1}$ is also compared with two kernel-free methods, namely SQSSVM and QSMPM.
	
	The penalty parameter $C$ for all methods
	is chosen from $\{10^{-7}, 10^{-6}, \ldots, 10^{7}\}$, the parameter $\sigma$ and the kernel parameter $p$ of all involved kernel functions 
	are chosen from $\{(\sqrt{2})^{-7}, (\sqrt{2})^{-6}, \ldots, (\sqrt{2})^{7}\}$, the parameter $\nu$ of $\nu$SVM is chosen from $\{0.1, 0.2, \ldots, 0.9\}$, the parameter $\tau_{1}$ of PSVM, the parameter $\mu$ of RSVM and the parameter $\kappa$ of RSHSVM are all selected from the set $\{0.1, 0.5, 0.9\}$. The 10-fold cross-validation technology is used to select the optimal parameters of each method. In addition, all features of each dataset is scaled to $[-1, 1]$ to ensure that the features have similar scales. The number of maximum iterations $K=10^{3}$, parameter $\eta=1.618$. Numerical experiments were run on a laptop with an Intel(R) Core(TM) i5-5200 CPU @ 2.20GHz and 8GB RAM.
	
	The P-stationary point is considered as a stopping criterion in the experiments, according to Theorem \ref{theorem}. We stop the algorithm when the point $(\widetilde{\boldsymbol{w}}^{k}; \boldsymbol{b}^{k}; c^{k}; \boldsymbol{u}^{k}; \boldsymbol{\lambda}^{k}))$ closely satisfies the conditions in formulas (\ref{l1})-(\ref{l5}), namely,
	\begin{equation}
		\max\{\theta_{1}^{k}, \theta_{2}^{k}, \theta_{3}^{k}, \theta_{4}^{k}\}<\operatorname{\tau},
	\end{equation} 
	where $\operatorname{\tau}$ represents the tolerance level, and the $\operatorname{\tau}=10^{-3}$ is set in the experiments. Specifically,
	\begin{equation}
		\begin{aligned}
			\theta_{1}^{k}
			&:=\frac{\|\boldsymbol{G}[\widetilde{\boldsymbol{w}}^{k}; \boldsymbol{b}^{k}]+[\boldsymbol{A}_{T_{k}}\:\boldsymbol{B}_{T_{k}}]^{\top}\boldsymbol{\lambda}_{T_{k}}^{k}\|}{1+\|[\widetilde{\boldsymbol{w}}^{k}; \boldsymbol{b}^{k}]\|},\\
			\theta_{2}^{k}
			&:=\frac{|\langle \boldsymbol{y}_{T_{k}}, \boldsymbol{\lambda}^{k}_{T_{k}}\rangle|}{1+|T_{k}|},\\
			\theta_{3}^{k}
			&:=\frac{\|\boldsymbol{u}^{k}-\boldsymbol{1}+\boldsymbol{A}\widetilde{\boldsymbol{w}}^{k}+\boldsymbol{B}\boldsymbol{b}^{k}+c^{k}\boldsymbol{y}\|}{\sqrt{N}},\\
			\theta_{4}^{k}
			&:=\frac{\|\boldsymbol{u}^{k}-\operatorname{prox}_{\gamma C\|(\cdot)_{+}\|_{0}}(\boldsymbol{u}^{k}-\gamma \boldsymbol{\lambda}^{k})\|}{1+\|\boldsymbol{u}^{k}\|}. \\
		\end{aligned}
	\end{equation}
	
	The numerical experiments are conducted on 4 artificial datasets and 14 benchmark datasets. For the experiments on artificial datasets, the visualization results are presented, the accuracy (ACC) and the number of SVs (NSV) are used to evaluate our QSSVM$_{0/1}$. For the experiments on benchmark datasets, the 10-fold cross-validation is repeated 10 times of each compared method, the mean value of ACC (mACC), the mean number of SVs (mNSV), and the CPU time are recorded.
	
	\subsection{Artificial datasets }
	In this subsection, the robustness, sample sparsity and effectiveness of our QSSVM$_{0/1}$ are explored. Four artificial datasets are constructed in two-dimensional space, which are Example 1, Example 2, Example 3 and Example 4, respectively. Each artificial dataset consists of 300 training samples with 2-dimensional features. During the processes of drawing plot, the bold red ``$\ast$'', blue ``$\circ$'', green ``$\circ$'', and bold black curve represent +1 class samples, -1 class samples, SVs, and hyperplanes or quadratic hyper-surfaces, respectively. The value of ACC or NSV is marked in each figure. Label noises or outliers are highlighted in magenta rectangular boxes.
	
		\begin{figure*}[!htbp]
		\centering
		\subfigure[Example 1]{
			\begin{minipage}[t]{0.2\linewidth}
				\centering
				\includegraphics[trim=36 0 36 0, clip, width=1.2in]{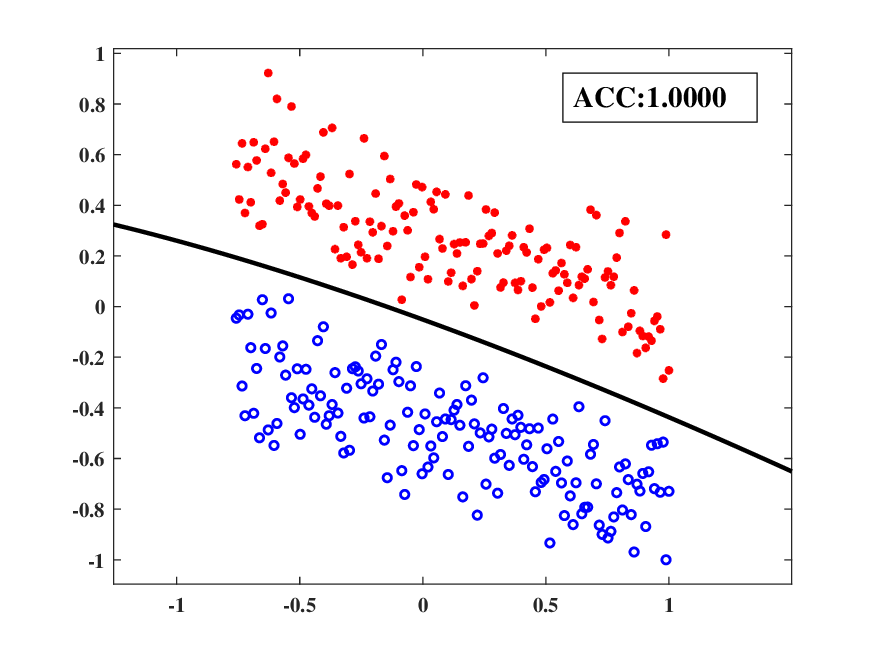}
			\end{minipage}
		}%
		\subfigure[Example 1 with 2 outliers]{
			\begin{minipage}[t]{0.2\linewidth}
				\centering
				\includegraphics[trim=36 0 36 0, clip, width=1.2in]{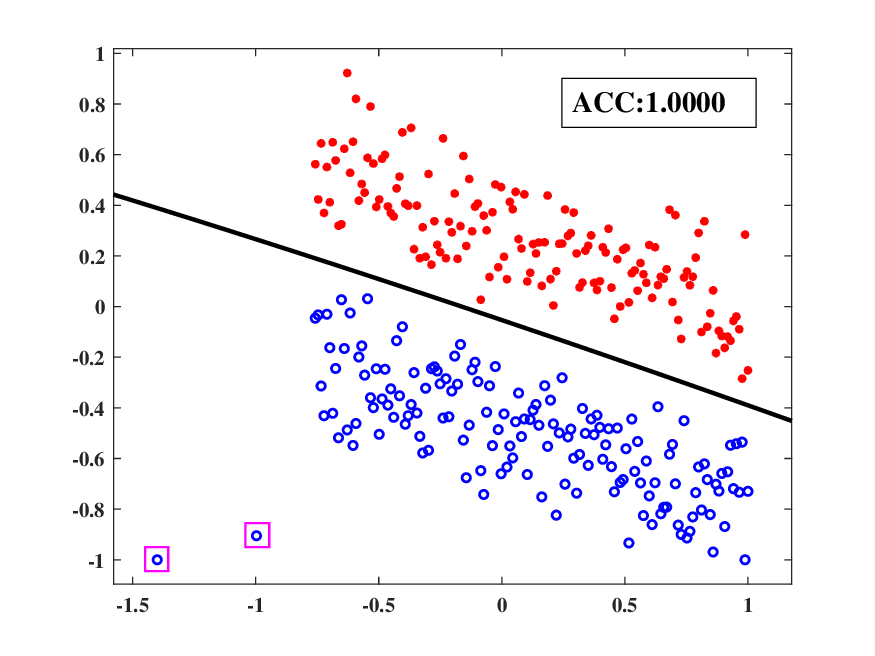}
			\end{minipage}
		}
		\subfigure[Example 1 with 2 label noises ]{
			\begin{minipage}[t]{0.2\linewidth}
				\centering
				\includegraphics[trim=36 0 36 0, clip, width=1.2in]{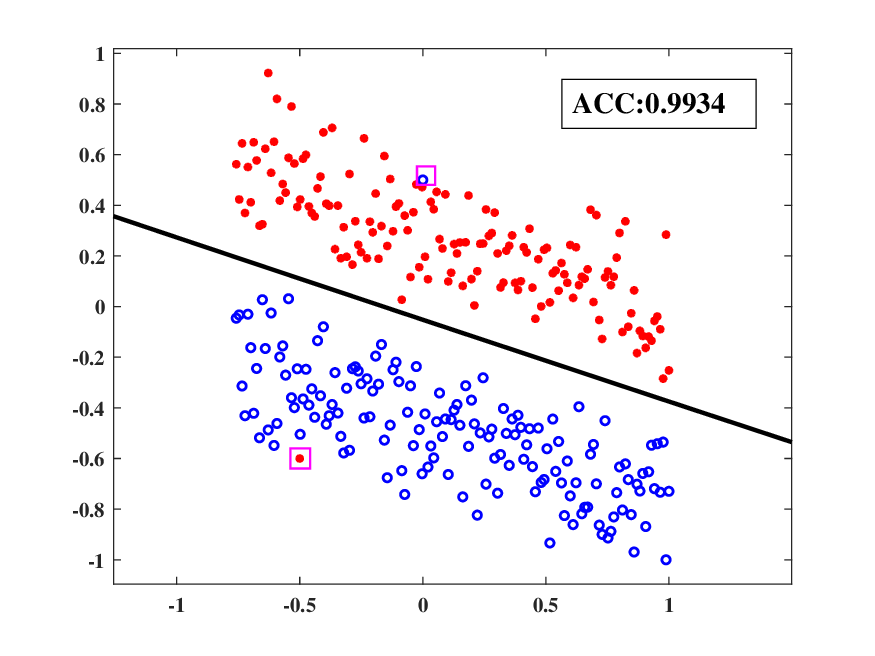}
			\end{minipage}
		}%
		\subfigure[Example 1 with outlier and label noise]{
			\begin{minipage}[t]{0.2\linewidth}
				\centering
				\includegraphics[trim=36 0 36 0, clip, width=1.2in]{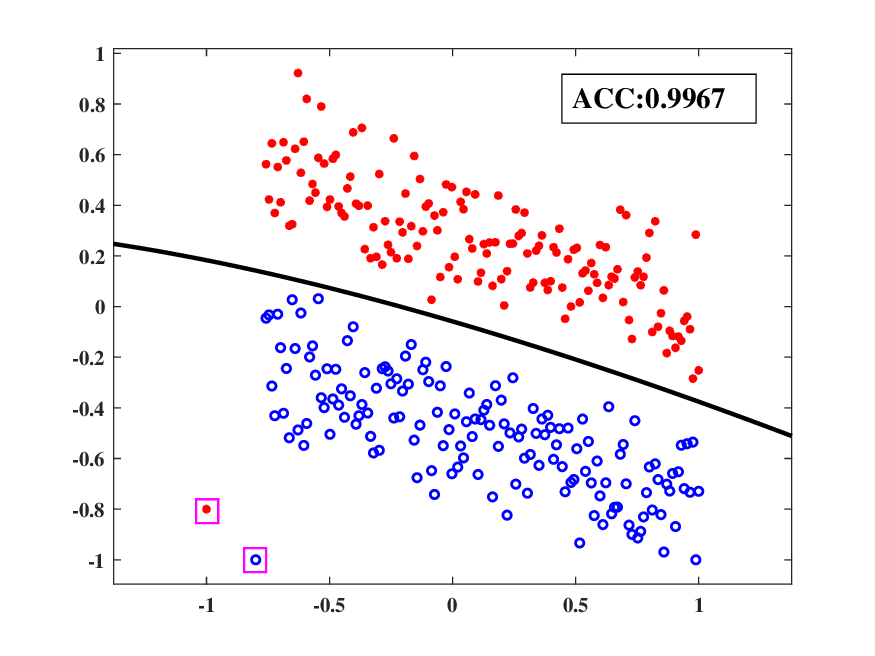}
			\end{minipage}
		}%
		\centering
		\caption{Classification results on Example 1 and Example 1 with label noises or outliers}\label{froubst}
	\end{figure*}
	
		\Cref{froubst} shows the visualization results of our QSSVM$_{0/1}$ on original Example 1 and the Example 1 with label noises or outliers. It can be observed that the QSSVM$_{0/1}$ is robust no matter adding label noises or some outliers.
		
\begin{figure*}[!t]
		\centering
		\subfigure[SVM$_{0/1}$]{
			\begin{minipage}[t]{0.2\linewidth}
				\centering
				\includegraphics[trim=36 0 36 0, clip, width=1.2in]{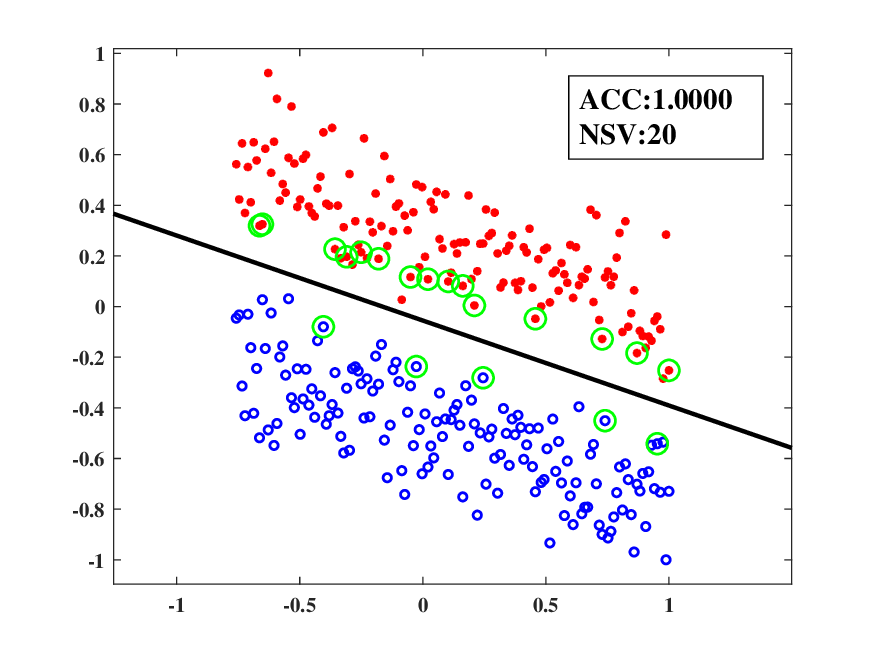}\label{figure a}
			\end{minipage}
		}%
		\subfigure[KSVM$_{0/1}$]{
			\begin{minipage}[t]{0.2\linewidth}
				\centering
				\includegraphics[trim=36 0 36 0, clip, width=1.2in]{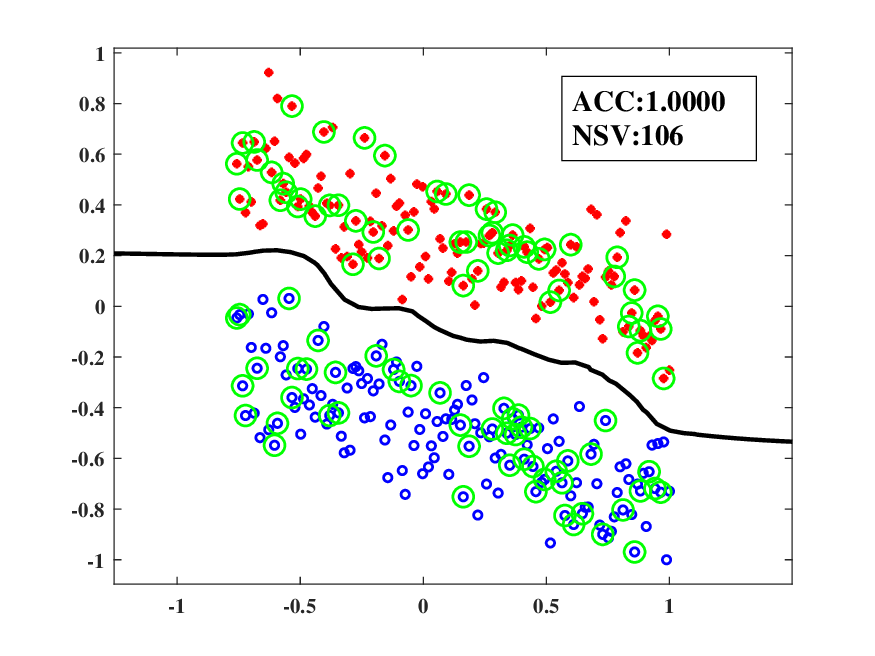}
			\end{minipage}
		}
		\subfigure[SQSSVM]{
			\begin{minipage}[t]{0.2\linewidth}
				\centering
				\includegraphics[trim=36 0 36 0, clip, width=1.2in]{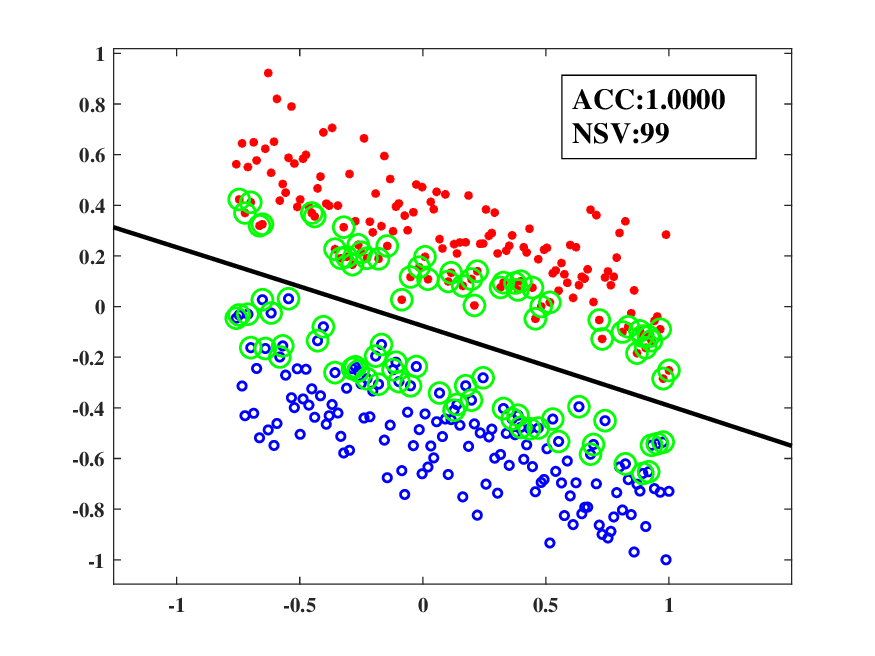}
			\end{minipage}
		}%
		\subfigure[QSSVM$_{0/1}$]{
			\begin{minipage}[t]{0.2\linewidth}
				\centering
				\includegraphics[trim=36 0 36 0,clip, width=1.2in]{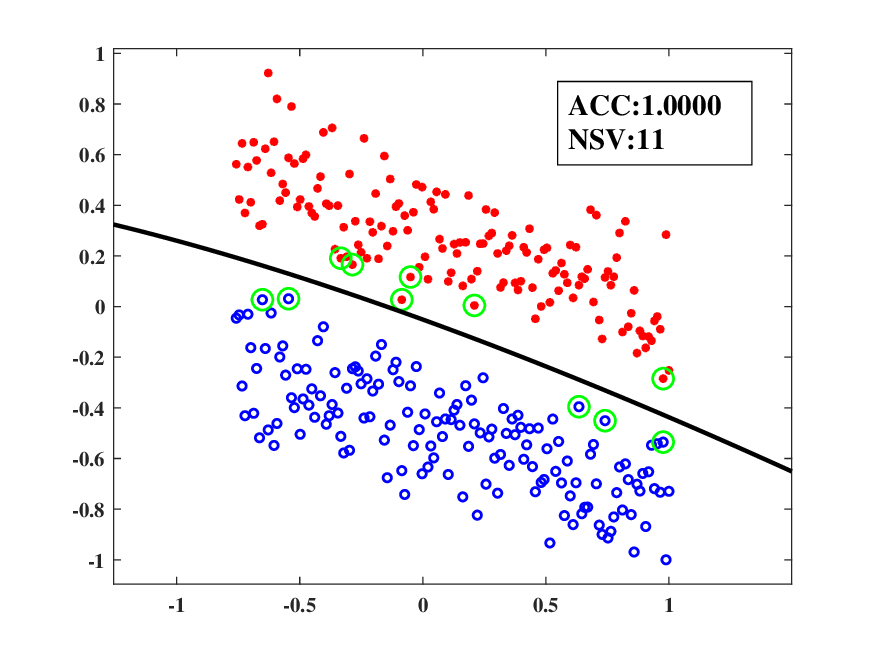}
			\end{minipage}
		}%
		\centering
		\caption{The classification results on Example 1}\label{f1}
	\end{figure*}

			The classification results of  SVM$_{0/1}$, KSVM$_{0/1}$, SQSSVM, and QSSVM$_{0/1}$ on the Example 1  are shown in \Cref{f1}. The values of ACC of all methods are 1, and the smooth quadratic hyper-surface obtained by QSSVM$_{0/1}$ approaches a straight line, but the hyper-surface obtained by KSVM$_{0/1}$ is not smooth. In addition, our QSSVM$_{0/1}$ has the lowest NSV.
		
\begin{figure*}[!htbp]
		\centering
		\subfigure[SVM$_{0/1}$]{
			\begin{minipage}[t]{0.2\linewidth}
				\centering
				\includegraphics[trim=36 0 36 0, clip, width=1.2in]{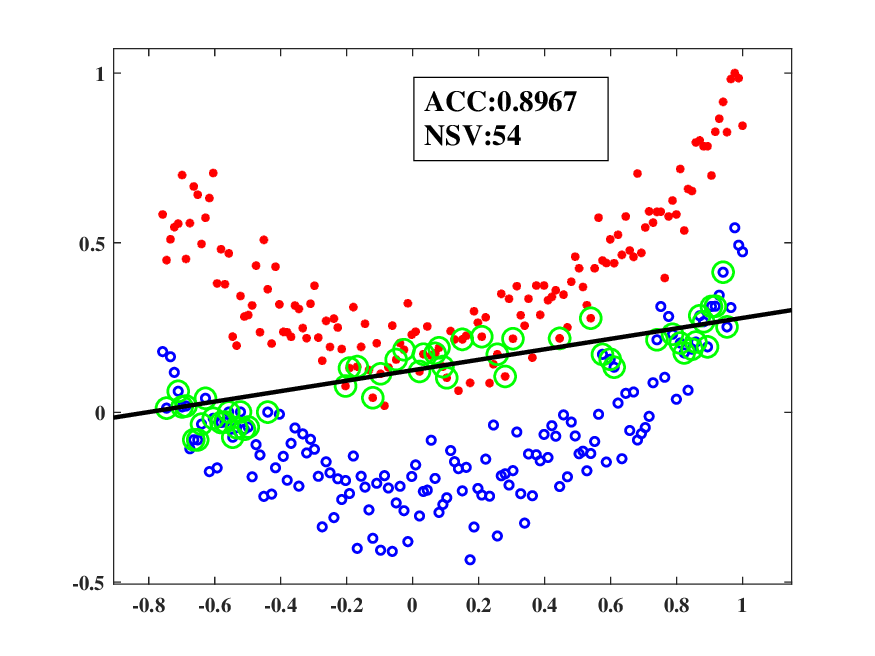}
			\end{minipage}
		}%
		\subfigure[KSVM$_{0/1}$]{
			\begin{minipage}[t]{0.2\linewidth}
				\centering
				\includegraphics[trim=36 0 36 0, clip, width=1.2in]{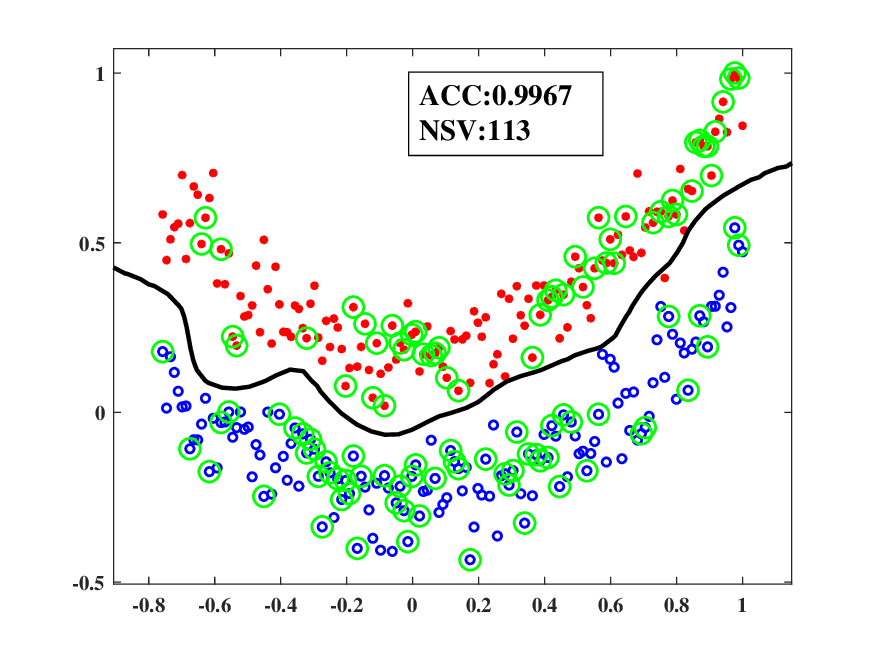}
			\end{minipage}
		}
		\subfigure[SQSSVM]{
			\begin{minipage}[t]{0.2\linewidth}
				\centering
				\includegraphics[trim=36 0 36 0, clip, width=1.2in]{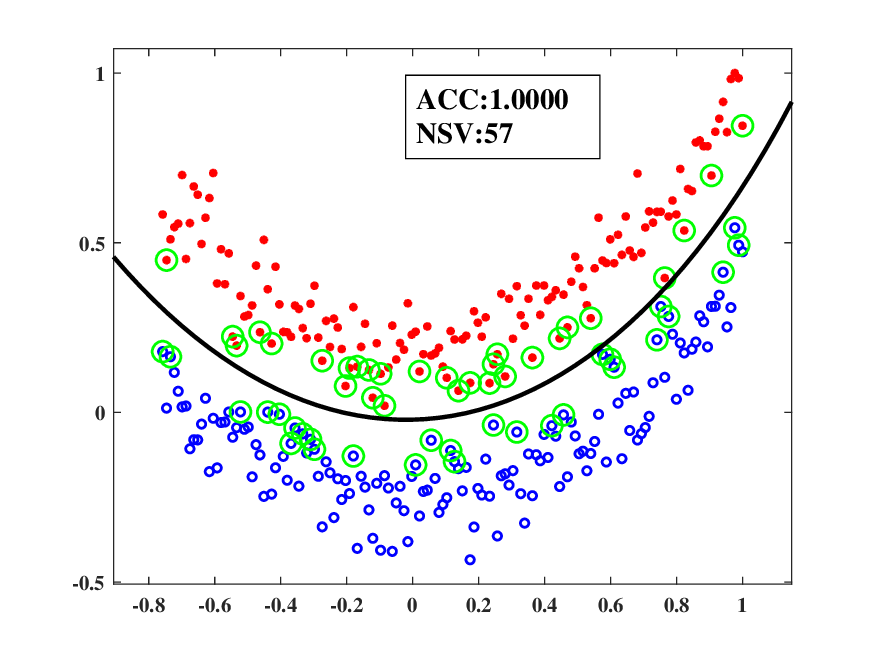}
			\end{minipage}
		}%
		\subfigure[QSSVM$_{0/1}$]{
			\begin{minipage}[t]{0.2\linewidth}
				\centering
				\includegraphics[trim=36 0 36 0, clip, width=1.2in]{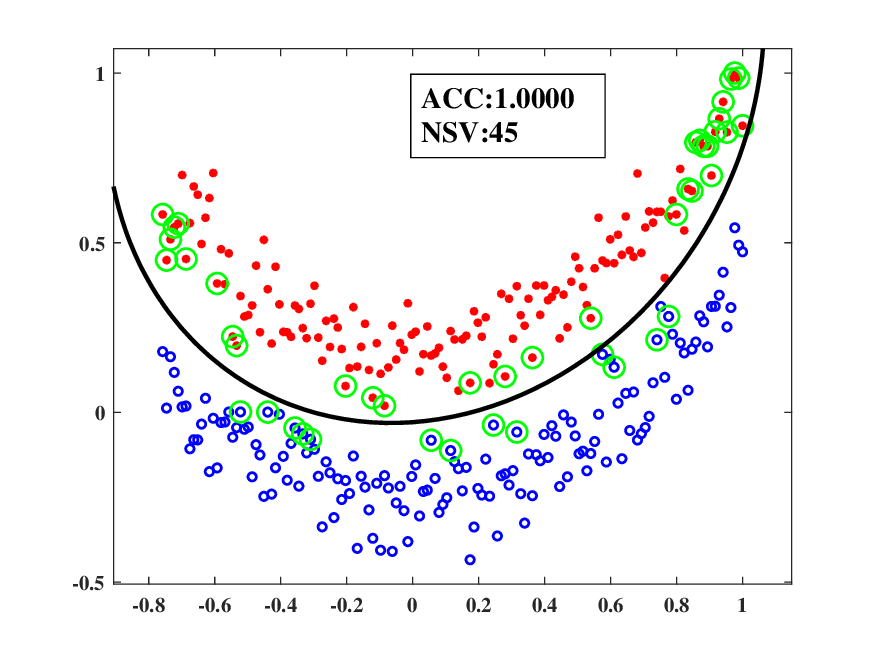}
			\end{minipage}
		}%
		\centering
		\caption{The classification results on Example 2}\label{f2}
	\end{figure*}

		\Cref{f2} reveals the classification results of the SVM$_{0/1}$, KSVM$_{0/1}$, SQSSVM, and QSSVM$_{0/1}$ on the Example 2. It is obvious that the ACC value of our QSSVM$_{0/1}$ is higher than that of SVM$_{0/1}$ and KSVM$_{0/1}$. The NSV of our QSSVM$_{0/1}$ is fewest. Furthermore, the smooth quadratic hyper-surface of QSSVM$_{0/1}$ is a parabola, but the hyper-surface of KSVM$_{0/1}$ is overfitting.
		
\begin{figure*}[!htbp]
		\centering
		\subfigure[SVM$_{0/1}$]{
			\begin{minipage}[t]{0.2\linewidth}
				\centering
				\includegraphics[trim=36 0 36 0, clip, width=1.2 in]{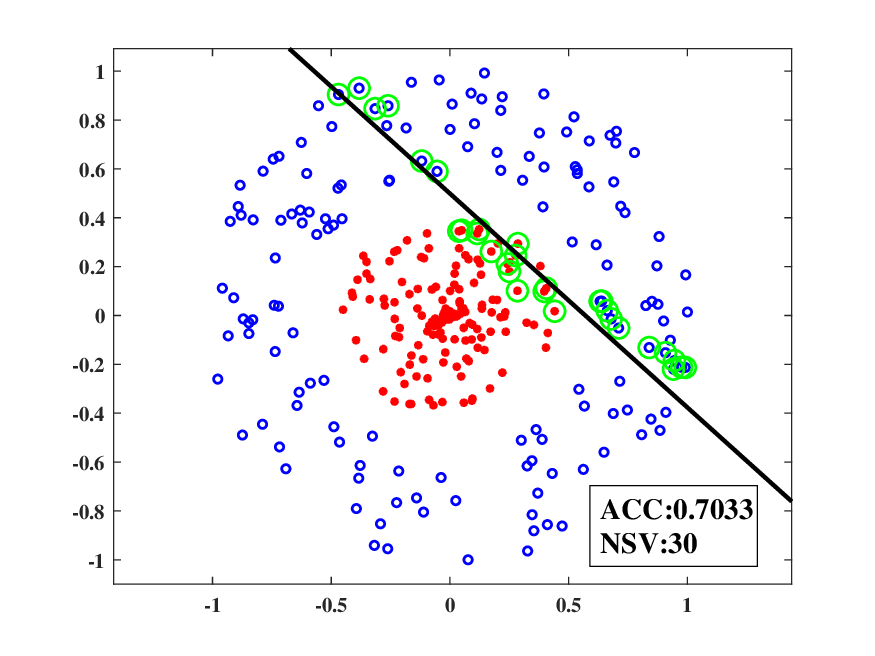}
			\end{minipage}
		}%
		\subfigure[KSVM$_{0/1}$]{
			\begin{minipage}[t]{0.2\linewidth}
				\centering
				\includegraphics[trim=36 0 36 0, clip, width=1.2in]{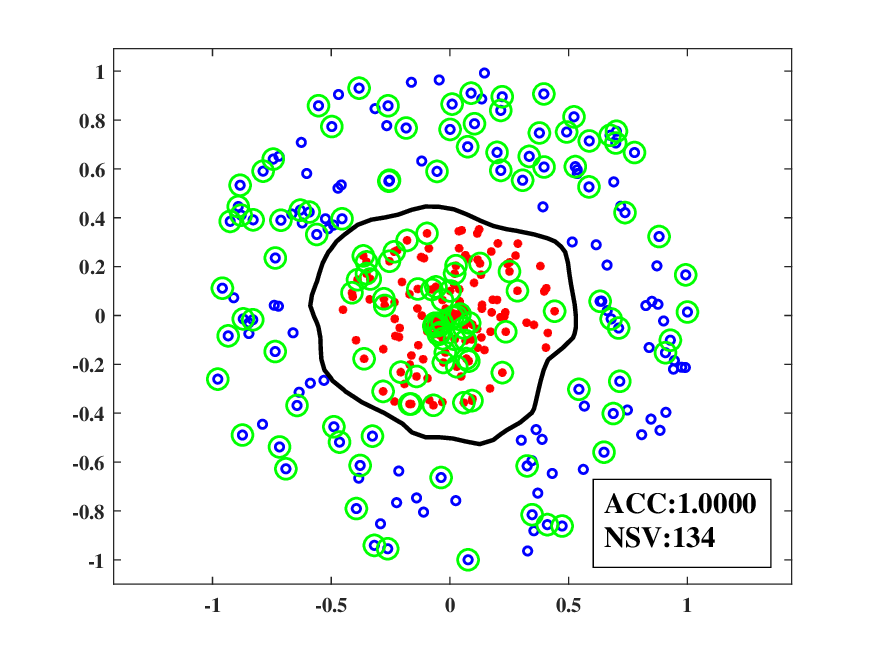}
			\end{minipage}
		}
		\subfigure[SQSSVM]{
			\begin{minipage}[t]{0.2\linewidth}
				\centering
				\includegraphics[trim=36 0 36 0, clip, width=1.2in]{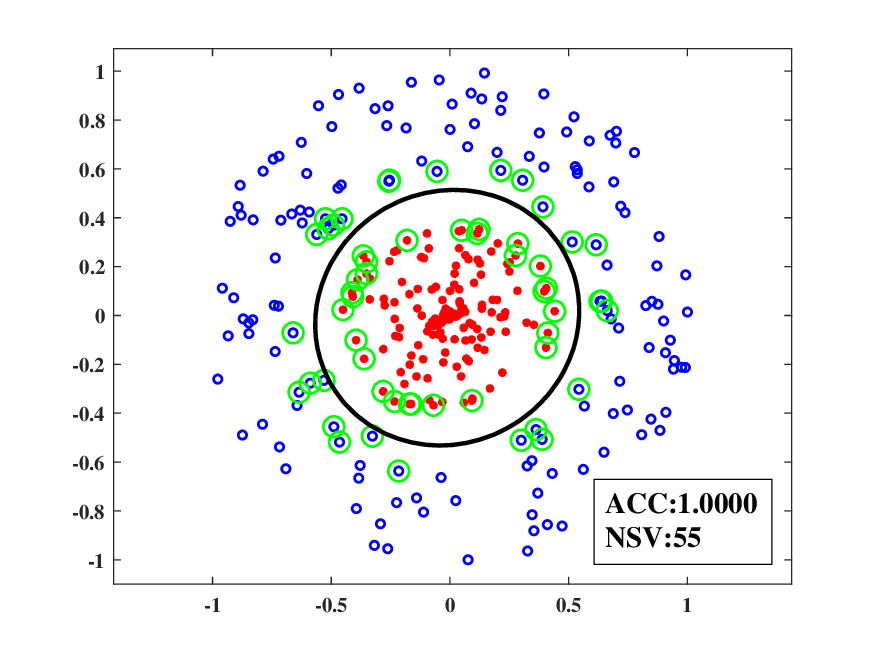}
			\end{minipage}
		}%
		\subfigure[QSSVM$_{0/1}$]{
			\begin{minipage}[t]{0.2\linewidth}
				\centering
				\includegraphics[trim=36 0 36 0, clip, width=1.2in]{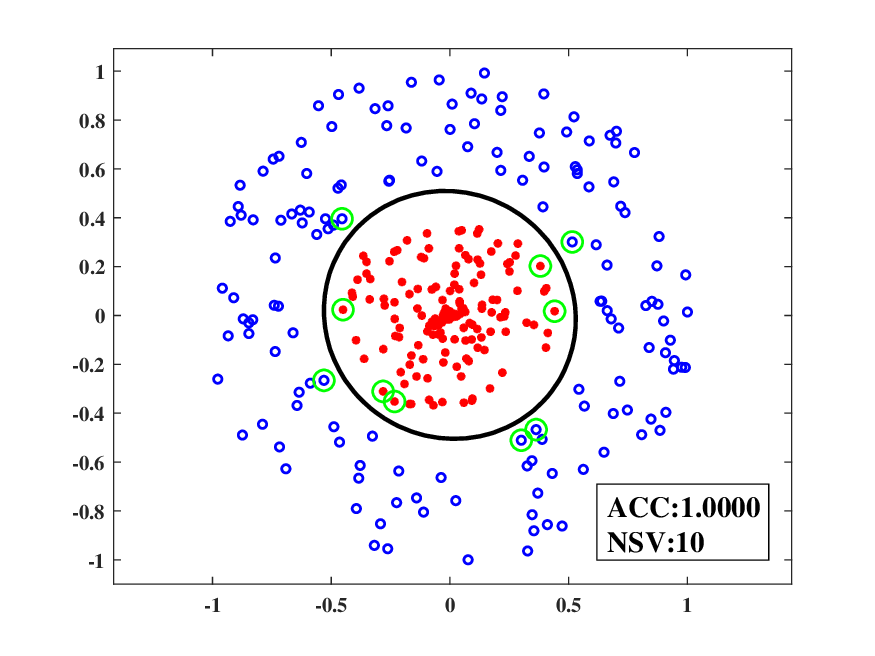}
			\end{minipage}
		}%
		\centering
		\caption{The classification results on Example 3}\label{f3}
	\end{figure*}

		The classification results of  SVM$_{0/1}$, KSVM$_{0/1}$, SQSSVM, and QSSVM$_{0/1}$ on the Example 3  are displayed in \Cref{f3}. It can be seen from \Cref{f3} that the ACC values of QSSVM$_{0/1}$, KSVM$_{0/1}$ and SQSSVM are all 1, and superior than that of SVM$_{0/1}$. Meanwhile, our QSSVM$_{0/1}$ has sample sparsity, because it has fewer SVs. Moreover, a circle as the separating hyper-surface is obtained by our QSSVM$_{0/1}$, and it is smooth.
		

\begin{figure*}[!t]
		\label{figure 5}
		\centering
		\subfigure[SVM$_{0/1}$]{
			\begin{minipage}[t]{0.2\linewidth}
				\centering
				\includegraphics[trim=36 0 36 0, clip, width=1.2in]{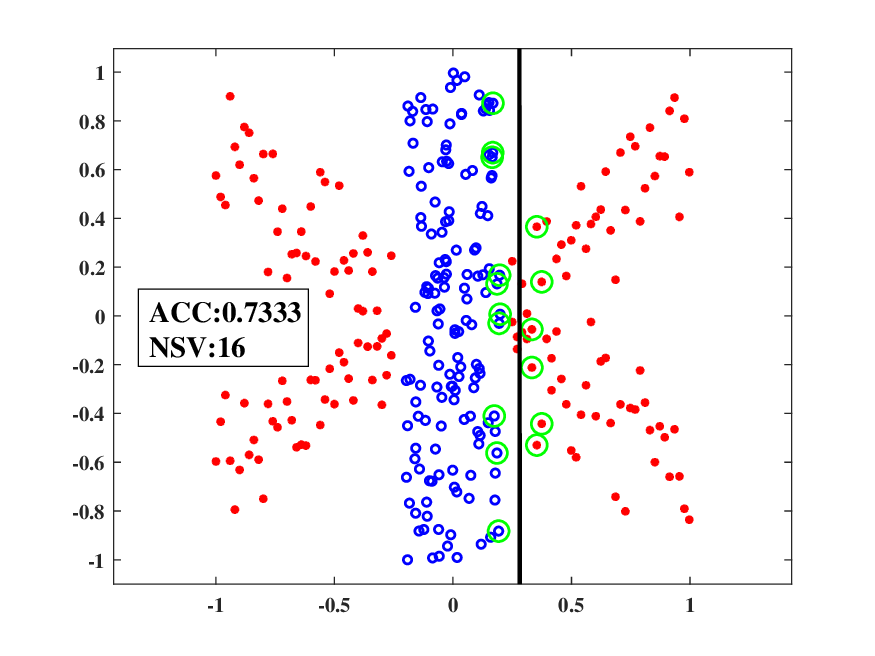}
			\end{minipage}
		}%
		\subfigure[KSVM$_{0/1}$]{
			\begin{minipage}[t]{0.2\linewidth}
				\centering
				\includegraphics[trim=36 0 36 0, clip, width=1.2in]{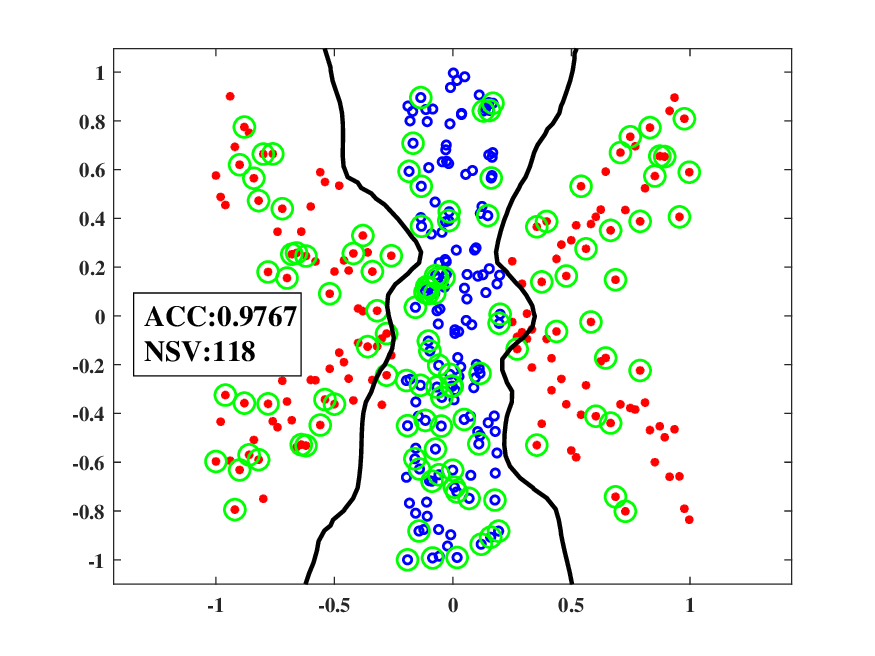}
			\end{minipage}
		}
		\subfigure[SQSSVM]{
			\begin{minipage}[t]{0.2\linewidth}
				\centering
				\includegraphics[trim=36 0 36 0, clip, width=1.2in]{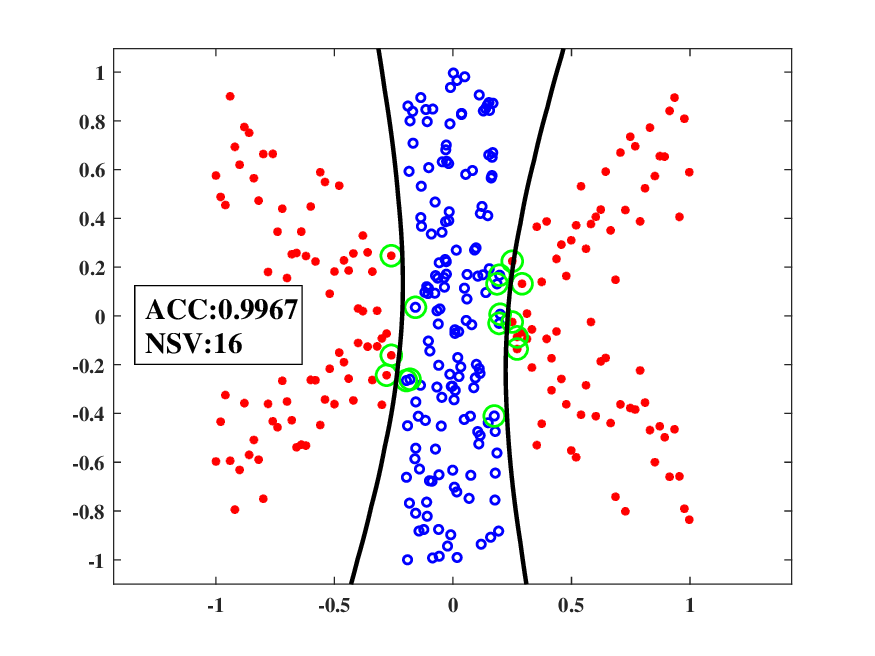}
			\end{minipage}
		}%
		\subfigure[QSSVM$_{0/1}$]{
			\begin{minipage}[t]{0.2\linewidth}
				\centering
				\includegraphics[trim=36 0 36 0, clip, width=1.2in]{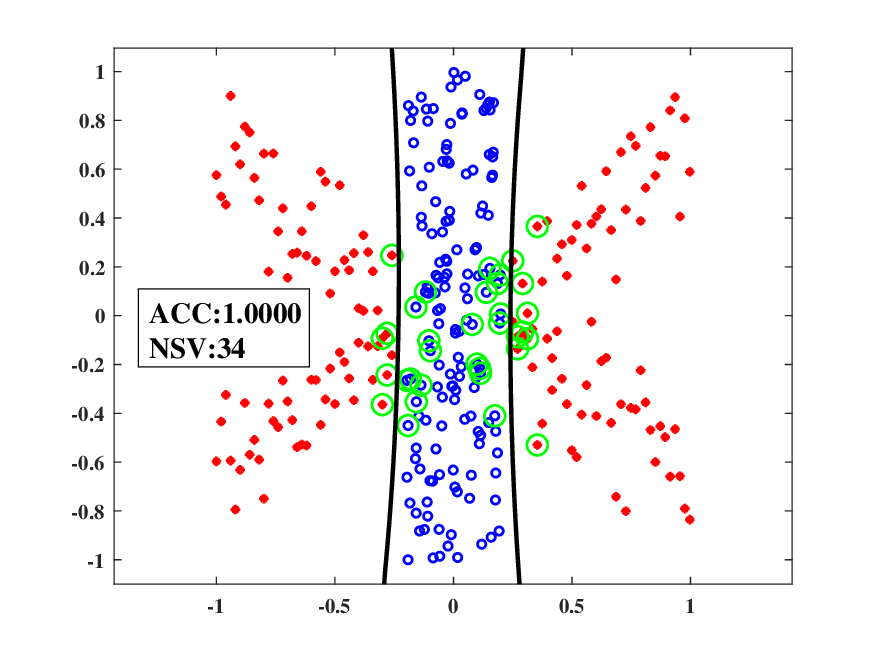}\label{figure t}
			\end{minipage}
		}%
		\centering
		\caption{The classification results on Example 4}\label{f5}
	\end{figure*}

	The classification results of SVM$_{0/1}$, KSVM$_{0/1}$, SQSSVM, and QSSVM$_{0/1}$ on the Example 4 are shown in \Cref{f5}. It can be seen that the value of ACC of our QSSVM$_{0/1}$ is better than that of the other methods. And QSSVM$_{0/1}$ has an advantage over KSVM$_{0/1}$ in terms of NSV. In addition, a hyperbola is obtained by QSSVM$_{0/1}$ to separate samples, but the hyper-surface of KSVM$_{0/1}$ is not smooth.

	Overall, a flexible separating hyper-surface can be obtained by our QSSVM$_{0/1}$, such as the line, parabola, circle and hyperbola displayed in the above numerical experiments. Furthermore it has higher values of ACC and fewer SVs, so it has better effectiveness and sample sparsity on artificial datasets.
	
		\subsection{Benchmark datasets}
	In addition to experiments on artificial datasets, QSSVM$_{0/1}$ is also compared with other methods on 14 benchmark datasets. The evaluation criteria mACC, mNSV and CPU time are used to verify the effectiveness, sample sparsity and efficiency of our method.  \Cref{table_1} provides a basic introduction to the 14 benchmark datasets. It is sorted based on the number of samples.
	
	\begin{table}[!htbp]
		\renewcommand{\arraystretch}{1}
		\caption{A basic information to 14 benchmark datasets}
		\centering
		\label{table_1}
		\small%
		\begin{tabular}{llllllll}
			\toprule
			\textbf{Datasets(Abbreviations)}&\textbf{Samples}&\textbf{Features}&\textbf{Classes}\\
			\midrule
			Caesarian(Cae) &80  &5 &2\\
			Tae(Tae)       &151 &5 &2\\
			Glass(Gla)     &214 &9 &2\\
			Heart-c(Hea) &303 &14&2\\
			Bupa(Bup)      &345 &6 &2\\
			Australian(Aus) &690 &14&2\\
			Pima(Pim)      &768 &8 &2\\
			Banknote(Ban)   &1372 &4 &2\\
			Yeast(Yea)      &1484 &8 &2\\
			Winequality(Win)&1599 &11&6\\
			Wireless(Wir)   &2000 &7 &4\\
			Abalone(Aba)    &2649 &8 &2\\
			Waveform(Wav)   &5000 &21&3\\
			Twonorm(Two)    &7400 &20&2\\
			\bottomrule
		\end{tabular}
	\end{table}
	
		\subsubsection{\bf Comparing with 0-1 loss function methods}
	
	On 14 benchmark datasets, our QSSVM$_{0/1}$ is compared to 2 methods with 0-1 loss function, which are SVM$_{0/1}$ and KSVM$_{0/1}$.
	The experimental results are represented in Table \ref{table_2}, and the optimal results are shown in bold.
	
	From Table \ref{table_2}, it can be concluded that the mACC of our QSSVM$_{0/1}$ outperforms that of the other two methods on the 12 benchmark datasets. Meanwhile its value of mNSV is fewer than that of the other methods on most datasets, such as on the Glass, Heart-c, Bupa, Pima, Abalone and Waveform datasets, respectively. As for CPU time, our QSSVM$_{0/1}$ costs less time than KSVM$_{0/1}$ on most datasets.
	
	\begin{table*}[!htbp]
		\tiny
		\centering
		\caption{The classification results of our QSSVM$_{0/1}$ and SVM$_{0/1}$, KSVM$_{0/1}$ }
	\label{table_2}
	\resizebox{\textwidth}{!}{
		\begin{tabular}{llllllllll}
			\toprule
			\multirow{2}{*}{\textbf{Datasets}} & \multicolumn{3}{c}{\textbf{mACC}(mean$\pm$std) } &\multicolumn{3}{c}{\textbf{mNSV}(mean$\pm$std)} &\multicolumn{3}{c}{\textbf{CPU time}(s)}  \\
			\cmidrule(lr){2-4}  \cmidrule(lr){5-7} \cmidrule(lr){8-10}
			&SVM$_{0/1}$&KSVM$_{0/1}$&QSSVM$_{0/1}$&SVM$_{0/1}$&KSVM$_{0/1}$&QSSVM$_{0/1}$&SVM$_{0/1}$&KSVM$_{0/1}$&QSSVM$_{0/1}$\\
			\midrule
			Cae&0.5925$\pm$0.0369&0.6425$\pm$0.0290&\textbf{0.6488}$\pm$0.0435&\textbf{2.6600}$\pm$0.5967&7.8800$\pm$0.2530&6.6300$\pm$0.4001
			&1.2667&\textbf{1.0312}&1.1981\\
			Tae      &0.7886$\pm$0.0186&0.7003$\pm$0.0779&\textbf{0.8140}$\pm$0.0113&10.8200$\pm$0.8364&\textbf{7.0800}$\pm$0.2530&9.5100$\pm$1.6031&1.2074&2.2761&\textbf{1.0374}\\
			Gla&0.6250$\pm$0.0274&0.6016$\pm$0.0360&\textbf{0.7255}$\pm$0.0354&15.2400$\pm$1.5806&18.5100$\pm$0.6488&\textbf{12.8200}$\pm$0.3327&\textbf{2.4617}&5.2416&6.8001\\
			Hea &0.9960$\pm$0.0115&0.9993$\pm$0.0021&\textbf{1.0000}$\pm$0.0000&23.5200$\pm$2.5459&30.2900$\pm$0.0316&\textbf{22.3000}$\pm$2.0499&\textbf{1.9765}&8.5582&23.2909\\
			Bup&0.6711$\pm$0.0126&0.5547$\pm$0.0293&\textbf{0.7043}$\pm$0.0135&23.5600$\pm$1.0875&13.9400$\pm$0.9879&\textbf{2.0000}$\pm$0.0000&\textbf{2.5974}&9.9653&6.2541\\
			Aus &\textbf{0.8584}$\pm$0.0060&0.8386$\pm$0.0122&0.8555$\pm$0.0027&28.4800$\pm$2.9253&\textbf{21.7500}$\pm$2.3524&35.4700$\pm$1.3630
			&\textbf{4.4616}&41.4978&31.1518\\
			Pim&\textbf{0.7625}$\pm$0.0050&0.6747$\pm$0.0262&0.7467$\pm$0.0195&36.1600$\pm$1.7469&15.4500$\pm$0.9914&\textbf{14.0000}$\pm$0.0000&\textbf{3.3384}&51.6956&15.1415\\
			Ban &0.9891$\pm$0.0018&0.9879$\pm$0.0038&\textbf{0.9929}$\pm$0.0023&\textbf{10.8200}$\pm$0.6250&108.9300$\pm$7.1824&44.9700$\pm$12.9220&\textbf{0.7925}&238.8360&6.7174\\
			Yea&0.6815$\pm$0.0093&0.6863$\pm$0.0085&\textbf{0.7073}$\pm$0.0072&100.6000$\pm$8.5723&\textbf{14.3600}$\pm$1.3729&117.0700$\pm$5.5148&\textbf{5.6410}&217.7618&21.2895\\
			Win&0.7226$\pm$0.0201&0.6834$\pm$0.0070&\textbf{0.7236}$\pm$0.0111&122.0600$\pm$5.0518&\textbf{29.3200}$\pm$2.1426&117.2800$\pm$5.8260&\textbf{9.0902}&267.3249&36.9722\\
			Wir&0.9908$\pm$0.0025&0.9682$\pm$0.0134&\textbf{0.9943}$\pm$0.0013&\textbf{15.9400}$\pm$2.1767&22.1700$\pm$1.4960&108.5500$\pm$13.9954&\textbf{1.8673}&369.8753&15.9246\\
			Aba&0.7623$\pm$0.0129&0.7221$\pm$0.0310&\textbf{0.8143}$\pm$0.0158&194.1700$\pm$16.9173&147.1700$\pm$5.4000&\textbf{55.4100}$\pm$15.8568&\textbf{8.6815}&1246.6851&14.8591\\
			Wav&0.8395$\pm$0.0105&0.8746$\pm$0.0019&\textbf{0.8874}$\pm$0.0049&148.5300$\pm$4.9878&499.7900$\pm$0.1287&\textbf{119.0000}$\pm$0.0000&\textbf{11.1868}&2663.9963&652.9437\\
			Two&0.9723$\pm$0.0137&0.9608$\pm$0.0022&\textbf{0.9732}$\pm$0.0007&99.1000$\pm$26.3912&\textbf{23.8500}$\pm$2.8289&227.5800$\pm$4.7719&\textbf{24.5967}&5550.9399&358.1596\\
			\bottomrule
	\end{tabular}}
	\end{table*}
	
		\subsubsection{\bf Comparing with other methods}
	
	Next, we compare our QSSVM$_{0/1}$ with SVM variants using different loss functions, including the SVM$_{L/R/P}$, $\nu$SVM$_{L/R/P}$, PSVM$_{L/R/P}$, LSSVM$_{L/R/P}$, RSVM$_{R}$ and RSHSVM$_{R}$. Here the subscripts $L, R, P$ represent linear, RBF and polynomial kernel functions, respectively. Additionally, the kernel-free classification methods SQSSVM and QSMPM are used to compare with our QSSVM$_{0/1}$.
	
	The numerical experimental results are displayed in \Cref{table_3}-\Cref{table_6}. The symbol ``$-$'' indicates that the classifier does not possess the corresponding evaluation criterion, while ``$--$'' denotes that the respective method fails to produce classification results due to computational time limitations or insufficient memory. The optimal results are highlighted in bold.
	
		\Cref{table_3} and \Cref{table_4} present the results of the comparison results between QSSVM$_{0/1}$ and 16 other methods with respect to the evaluation criterion mACC. From the results of these two tables, it can be concluded that our QSSVM$_{0/1}$ has higher accuracy than that of other methods on most datasets. This demonstrates its effectiveness in solving classification problems.
	
	The results of the 17 methods with respect to the mNSV evaluation criterion are illustrated in \Cref{table_5} and \Cref{table_6}. By comparison, it is found that our QSSVM$_{0/1}$ achieves the minimum value in terms of the mNSV on 6 datasets, which verifies the sample sparsity of our QSSVM$_{0/1}$.
		
	\begin{table*}[!htbp]
			\tiny
		\centering
		\caption{mACC classification results (mean$\pm$std) on 7 small-scale datasets}
		\label{table_3}
		\resizebox{\textwidth}{!}{
			\begin{tabular}{llllllllll}
				\toprule
				\textbf{Models}&\textbf{Cae}&\textbf{Tae}&\textbf{Gla}&\textbf{Hea}
				&\textbf{Bup}&\textbf{Aus}&\textbf{Pim} \\
				\midrule
				SVM$_{L}$&0.6250$\pm$0.0156&0.8031$\pm$0.0104&0.6590$\pm$0.0114&\textbf{1.0000}$\pm$0.0000&0.6847$\pm$0.0124&0.8551$\pm$0.0000&0.6511$\pm$0.0001\\
				SVM$_{R}$&0.6438$\pm$0.0259&0.8015$\pm$0.0132&0.7619$\pm$0.0174&\textbf{1.0000}$\pm$0.0000&0.7022$\pm$0.0106&0.8546$\pm$0.0031&0.7436$\pm$0.0040\\
				SVM$_{P}$&0.6388$\pm$0.0375&0.8094$\pm$0.0161&0.7487$\pm$0.0166&\textbf{1.0000}$\pm$0.0000&\textbf{0.7279}$\pm$0.0085&0.8658$\pm$0.0027&0.7720$\pm$0.0034\\
				$\nu$SVM$_{L}$&0.5750$\pm$0.0000&0.8071$\pm$0.0000&0.6449$\pm$0.0006&0.7976$\pm$0.0099&0.5797$\pm$0.0004&0.8604$\pm$0.0045&0.6510$\pm$0.0000\\
				$\nu$SVM$_{R}$&0.6000$\pm$0.0300&0.8104$\pm$0.0151&\textbf{0.8146}$\pm$0.0121&\textbf{1.0000}$\pm$0.0000&0.7001$\pm$0.0197&0.8633$\pm$0.0042&0.7661$\pm$0.0062\\
				$\nu$SVM$_{P}$&0.5950$\pm$0.0340&\textbf{0.8140}$\pm$0.0123&0.7490$\pm$0.0087&\textbf{1.0000}$\pm$0.0000&0.7272$\pm$0.0070&0.8568$\pm$0.0015&\textbf{0.7773}$\pm$0.0026\\
				PSVM$_{L}$&0.5988$\pm$0.0309&0.8082$\pm$0.0007&0.6466$\pm$0.0143&\textbf{1.0000}$\pm$0.0000&0.6771$\pm$0.0175&0.8551$\pm$0.0000&0.7628$\pm$0.0045\\
				PSVM$_{R}$&0.6325$\pm$0.0237&\textbf{0.8140}$\pm$0.0197&0.8000$\pm$0.0109&\textbf{1.0000}$\pm$0.0000&0.7174$\pm$0.0057&0.8586$\pm$0.0042&0.7750$\pm$0.0042\\
				PSVM$_{P}$&0.6250$\pm$0.0212&0.8099$\pm$0.0149&0.7434$\pm$0.0162&\textbf{1.0000}$\pm$0.0000&0.7189$\pm$0.0088&0.8662$\pm$0.0037&0.7763$\pm$0.0028\\
				LSSVM$_{L}$&0.6375$\pm$0.0228&0.8118$\pm$0.0098&0.6329$\pm$0.0116&\textbf{1.0000}$\pm$0.0000&0.6853$\pm$0.0072&\textbf{0.8668}$\pm$0.0043&0.6512$\pm$0.0008\\
				LSSVM$_{R}$&0.5525$\pm$0.0407&0.8139$\pm$0.0100&0.7613$\pm$0.0152&0.8174$\pm$000.71&0.7204$\pm$0.0090&0.8559$\pm$0.0017&0.7690$\pm$0.0038\\
				LSSVM$_{P}$&0.5988$\pm$0.0314&0.8130$\pm$0.0120&0.7668$\pm$0.0124&\textbf{1.0000}$\pm$0.0000&0.7262$\pm$0.0085&0.8648$\pm$0.0033&0.7688$\pm$0.0057\\
				RSVM$_{R}$&0.6188$\pm$0.0340&0.8113$\pm$0.0175&0.7240$\pm$0.0160&\textbf{1.0000}$\pm$0.0000&0.7003$\pm$0.0139&0.8551$\pm$0.0051&0.7531$\pm$0.0055\\
				RSHSVM$_{R}$&0.5750$\pm$0.0000&0.8076$\pm$0.0005&0.7319$\pm$0.0103&\textbf{1.0000}$\pm$0.0000&0.6961$\pm$0.0104&0.8523$\pm$0.0042&0.7557$\pm$0.0005\\
				SQSSVM&0.5925$\pm$0.0369&0.8079$\pm$0.0006&0.7324$\pm$0.0132&\textbf{1.0000}$\pm$0.0000&0.7206$\pm$0.0110&0.8545$\pm$0.0017&0.7663$\pm$0.0037\\
				QSMPM&0.6025$\pm$0.0293&0.6799$\pm$0.0189&0.7326$\pm$0.0123&\textbf{1.0000}$\pm$0.0000&0.7124$\pm$0.0045&0.8565$\pm$0.0047&0.7424$\pm$0.0062\\
				QSSVM$_{0/1}$&\textbf{0.6488}$\pm$0.0435&\textbf{0.8140}$\pm$0.0113&0.7255$\pm$0.0354&\textbf{1.0000}$\pm$0.0000&0.7043$\pm$0.0135&0.8555$\pm$0.0027&0.7467$\pm$0.0195\\
				\bottomrule
		\end{tabular}}
	\end{table*}
	
	\begin{table*}[!htbp]
			\tiny
		\centering
		\caption{mACC classification results (mean$\pm$std) on 7 large-scale datasets}
		\label{table_4}
		\resizebox{\textwidth}{!}{
			\begin{tabular}{llllllllll}
				\toprule
				\textbf{Models}&\textbf{Ban}&\textbf{Yea}&\textbf{Win}&\textbf{Wir}&\textbf{Aba}&\textbf{Wav}
				&\textbf{Two}\\
				\midrule
				SVM$_{L}$&0.9886$\pm$0.0007&0.6880$\pm$0.0000&0.7396$\pm$0.0028&0.9589$\pm$0.0302&0.8015$\pm$0.0007&0.8581$\pm$0.0010&\textbf{0.9781}$\pm$0.0003\\
				SVM$_{R}$&\textbf{1.0000}$\pm$0.0000&0.6983$\pm$0.0052&\textbf{0.7817}$\pm$0.0053&0.9901$\pm$0.0004&0.8264$\pm$0.0039&0.8832$\pm$0.0025&$--$\\
				SVM$_{P}$&\textbf{1.0000}$\pm$0.0000&0.7213$\pm$0.0055&0.7469$\pm$0.0030&0.9926$\pm$0.0003&0.8299$\pm$0.0016&0.8951$\pm$0.0008&$--$\\
				$\nu$SVM$_{L}$&0.6146$\pm$0.0601&0.6880$\pm$0.0000&0.5347$\pm$0.0000&0.7500$\pm$0.0000&0.5280$\pm$0.0303&0.6686$\pm$0.0000&0.8432$\pm$0.0469\\
				$\nu$SVM$_{R}$&0.9856$\pm$0.0017&0.7359$\pm$0.0044&0.7461$\pm$0.0038&0.9900$\pm$0.0003&\textbf{0.8376}$\pm$0.0021&0.8923$\pm$0.0010&$--$\\
				$\nu$SVM$_{P}$&0.9956$\pm$0.0000&0.7301$\pm$0.0049&0.7516$\pm$0.0028&0.9941$\pm$0.0004&0.8272$\pm$0.0015&0.8918$\pm$0.0006&$--$\\
				PSVM$_{L}$&0.9768$\pm$0.0017&0.6880$\pm$0.0000&0.6555$\pm$0.0058&0.9904$\pm$0.0005&0.8127$\pm$0.0016&0.8487$\pm$0.0070&0.9779$\pm$0.0003\\
				PSVM$_{R}$&0.9993$\pm$0.0000&\textbf{0.7373}$\pm$0.0042&0.7748$\pm$0.0038&0.9910$\pm$0.0002&0.8374$\pm$0.0023&0.8882$\pm$0.0013&$--$\\
				PSVM$_{P}$&\textbf{1.0000}$\pm$0.0000&0.6981$\pm$0.0055&0.7464$\pm$0.0033&0.9925$\pm$0.0004&0.8285$\pm$0.0009&0.8935$\pm$0.0013&$--$\\
				LSSVM$_{L}$&0.9765$\pm$0.0005&0.6761$\pm$0.0030&0.7417$\pm$0.0022&0.9858$\pm$0.0006&0.8015$\pm$0.0014&0.8562$\pm$0.0013&0.9780$\pm$0.0004
				\\
				LSSVM$_{R}$&0.9893$\pm$0.0004&0.6898$\pm$0.0045&0.7539$\pm$0.0021&0.9866$\pm$0.0002&0.8367$\pm$0.0012&0.8958$\pm$0.0010&$--$\\
				LSSVM$_{P}$&\textbf{1.0000}$\pm$0.0000&0.7132$\pm$0.0027&0.7472$\pm$0.0025&0.9937$\pm$0.0003&0.8301$\pm$0.0015&\textbf{0.9001}$\pm$0.0015&$--$\\
				RSVM$_{R}$&\textbf{1.0000}$\pm$0.0000&0.7057$\pm$0.0029&0.7542$\pm$0.0025&0.9873$\pm$0.0046&0.8009$\pm$0.0016&$--$&$--$\\
				RSHSVM$_{R}$&\textbf{1.0000}$\pm$0.0000&0.6880$\pm$0.0000&0.7123$\pm$0.0022&0.9931$\pm$0.0005&0.8308$\pm$0.0024&$--$&$--$\\
				SQSSVM&0.9847$\pm$0.0001&0.7206$\pm$0.0020&0.7414$\pm$0.0035&0.9934$\pm$0.0002&0.8292$\pm$0.0014&0.8842$\pm$0.0023&$--$\\
				QSMPM&0.9906$\pm$0.0013&0.6969$\pm$0.0029&0.7472$\pm$0.0019&0.9941$\pm$0.0006&0.8301$\pm$0.0010&0.8826$\pm$0.0013&0.9758$\pm$0.0003\\
				QSSVM$_{0/1}$&0.9929$\pm$0.0023&0.7073$\pm$0.0072&0.7236$\pm$0.0111&\textbf{0.9943}$\pm$0.0013&0.8143$\pm$0.0158&0.8874$\pm$0.0049&0.9732$\pm$0.0007\\
				\bottomrule
		\end{tabular}}
	\end{table*}
	
	\begin{table*}[!htbp]
		\tiny
		\centering
		\caption{mNSV classification results (mean$\pm$std) on 7 small-scale datasets}
		\label{table_5}
		\resizebox{\textwidth}{!}{
			\begin{tabular}{llllllllll}
				\toprule
				\textbf{Models}&\textbf{Cae}&\textbf{Tae}&\textbf{Gla}&\textbf{Hea}
				&\textbf{Bup}&\textbf{Aus}&\textbf{Pim} \\
				\midrule
				SVM$_{L}$&\textbf{5.7900}$\pm$0.5724&7.7300$\pm$0.7119&\textbf{10.2700}$\pm$0.5229&22.3900$\pm$0.6297&22.4700$\pm$1.0573&61.6700$\pm$0.7056&66.3400$\pm$2.1614\\
				SVM$_{R}$&8.0000$\pm$0.0000&14.8700$\pm$0.1252&18.1500$\pm$0.4552&30.3000$\pm$0.0000&30.7200$\pm$0.5750&69.0000$\pm$0.0000&65.1100$\pm$0.3071\\
				SVM$_{P}$&6.9100$\pm$0.2558&\textbf{6.7600}$\pm$0.3438&10.7700$\pm$0.8274&22.9800$\pm$0.3084&20.5600$\pm$0.9131&69.0000$\pm$0.0000&51.5000$\pm$0.6307\\
				$\nu$SVM$_{L}$&6.0700$\pm$0.5478&10.9000$\pm$0.0000&21.2300$\pm$0.3622&30.2700$\pm$0.0483&34.5000$\pm$0.0000&66.9900$\pm$0.4433&76.8000$\pm$0.0000\\
				$\nu$SVM$_{R}$&7.7000$\pm$0.1247&14.9900$\pm$0.0738&19.1600$\pm$0.3658&30.1800$\pm$0.1135&34.5000$\pm$0.0000&68.9900$\pm$0.0316&76.8000$\pm$0.0000\\
				$\nu$SVM$_{P}$&5.9100$\pm$0.3985&7.9800$\pm$0.4826&16.5900$\pm$0.8987&22.3200$\pm$0.4780&27.7800$\pm$0.3706&57.9900$\pm$0.1595&54.5900$\pm$1.8107\\
				PSVM$_{L}$&8.0000$\pm$0.0000&15.1000$\pm$0.0000&21.4000$\pm$0.0000&30.3000$\pm$0.0000&34.5000$\pm$0.0000&69.0000$\pm$0.0000&76.8000$\pm$0.0000\\
				PSVM$_{R}$&8.0000$\pm$0.0000&15.1000$\pm$0.0000&21.4000$\pm$0.0000&30.3000$\pm$0.0000&34.5000$\pm$0.0000&69.0000$\pm$0.0000&76.8000$\pm$0.0000\\
				PSVM$_{P}$&8.0000$\pm$0.0000&15.1000$\pm$0.0000&21.4000$\pm$0.0000&30.3000$\pm$0.0000&34.5000$\pm$0.0000&69.0000$\pm$0.0000&76.8000$\pm$0.0000\\
				LSSVM$_{L}$&8.0000$\pm$0.0000&15.1000$\pm$0.0000&21.4000$\pm$0.0000&30.3000$\pm$0.0000&34.5000$\pm$0.0000&69.0000$\pm$0.0000&76.8000$\pm$0.0000\\
				LSSVM$_{R}$&8.0000$\pm$0.0000&15.1000$\pm$0.0000&21.4000$\pm$0.0000&30.3000$\pm$0.0000&34.5000$\pm$0.0000&69.0000$\pm$0.0000&76.8000$\pm$0.0000\\
				LSSVM$_{P}$&8.0000$\pm$0.0000&15.1000$\pm$0.0000&21.4000$\pm$0.0000&30.3000$\pm$0.0000&34.5000$\pm$0.0000&69.0000$\pm$0.0000&76.8000$\pm$0.0000\\
				RSVM$_{R}$&5.8600$\pm$1.0585&11.5000$\pm$0.6912&11.1500$\pm$0.4882&29.9900$\pm$0.1287&18.3600$\pm$0.8922&45.4200$\pm$1.1679&14.9500$\pm$1.6847\\
				RSHSVM$_{R}$&7.8400$\pm$0.3373&12.7100$\pm$0.9445&15.5400$\pm$0.4427&25.2300$\pm$0.6075&30.1300$\pm$0.6395&44.6500$\pm$1.0416&62.1000$\pm$1.0934\\
				SQSSVM&6.1700$\pm$0.2830&8.8100$\pm$0.4149&11.1600$\pm$0.5275&26.3600$\pm$0.5661&19.2200$\pm$0.5846&42.0200$\pm$0.7772&44.9300$\pm$0.8642\\
				QSMPM&$-$&$-$&$-$&$-$&$-$&$-$&$-$\\
				QSSVM$_{0/1}$&6.6300$\pm$0.4001&9.5100$\pm$1.6031&12.8200$\pm$0.3327&\textbf{22.3000}$\pm$2.0499&\textbf{2.0000}$\pm$0.0000&\textbf{35.4700}$\pm$1.3630&\textbf{14.0000}$\pm$0.0000\\
				\bottomrule
		\end{tabular}}
	\end{table*}
	
	\begin{table*}[!htbp]
		\tiny
		\centering
		\caption{mNSV classification results (mean$\pm$std) on 7 large-scale datasets}
		\label{table_6}
		\resizebox{\textwidth}{!}{
			\begin{tabular}{llllllllll}
				\toprule
				\textbf{Models}&\textbf{Ban}&\textbf{Yea}&\textbf{Win}&\textbf{Wir}&\textbf{Aba}&\textbf{Wav}
				&\textbf{Two}\\
				\midrule
				SVM$_{L}$&8.8500$\pm$0.1900&148.4000$\pm$0.0000&91.7900$\pm$0.8863&101.4200$\pm$0.2781&254.6600$\pm$2.1578&204.7700$\pm$1.0478&\textbf{38.5200}$\pm$0.9624\\
				SVM$_{R}$&42.5200$\pm$0.4662&77.8400$\pm$1.5248&154.0000$\pm$0.4497&111.5200$\pm$0.3824&85.8700$\pm$1.9420&136.7500$\pm$1.6795&$--$\\
				SVM$_{P}$&14.9700$\pm$0.3743&87.2100$\pm$1.1986&86.1200$\pm$1.0250&105.4600$\pm$0.5190&106.8000$\pm$0.8179&308.3100$\pm$1.4255&$--$\\
				$\nu$SVM$_{L}$&65.3700$\pm$5.8867&148.4000$\pm$0.0000&159.9000$\pm$0.0000&194.3400$\pm$6.4032&242.3100$\pm$2.7950&500.0000$\pm$0.0000&668.3600$\pm$0.2171\\
				$\nu$SVM$_{R}$&81.6600$\pm$1.6249&131.9800$\pm$1.0675&159.9000$\pm$0.0000&188.3400$\pm$7.9200&264.9000$\pm$0.0000&448.1600$\pm$6.9666&$--$\\
				$\nu$SVM$_{P}$&134.6900$\pm$4.2800&148.4000$\pm$0.0000&150.2400$\pm$4.8484&69.8400$\pm$4.1000&264.9000$\pm$0.0000&263.0600$\pm$1.0113&$--$\\
				PSVM$_{L}$&137.2000$\pm$0.0000&148.4000$\pm$0.0000&159.9000$\pm$0.0000&200.0000$\pm$0.0000&264.9000$\pm$0.0000&500.0000$\pm$0.0000&740.0000$\pm$0.0000\\
				PSVM$_{R}$&137.2000$\pm$0.0000&148.4000$\pm$0.0000&159.9000$\pm$0.0000&200.0000$\pm$0.0000&264.9000$\pm$0.0000&500.0000$\pm$0.0000&740.0000$\pm$0.0000\\
				PSVM$_{P}$&137.2000$\pm$0.0000&148.4000$\pm$0.0000&159.9000$\pm$0.0000&200.0000$\pm$0.0000&264.9000$\pm$0.0000&500.0000$\pm$0.0000&740.0000$\pm$0.0000\\
				LSSVM$_{L}$&137.2000$\pm$0.0000&148.4000$\pm$0.0000&159.9000$\pm$0.0000&200.0000$\pm$0.0000&264.9000$\pm$0.0000&500.0000$\pm$0.0000&740.0000$\pm$0.0000\\
				LSSVM$_{R}$&137.2000$\pm$0.0000&148.4000$\pm$0.0000&159.9000$\pm$0.0000&200.0000$\pm$0.0000&264.9000$\pm$0.0000&500.0000$\pm$0.0000&740.0000$\pm$0.0000\\
				LSSVM$_{P}$&137.2000$\pm$0.0000&148.4000$\pm$0.0000&159.9000$\pm$0.0000&200.0000$\pm$0.0000&264.9000$\pm$0.0000&500.0000$\pm$0.0000&740.0000$\pm$0.0000\\
				RSVM$_{R}$&22.0800$\pm$0.7315&\textbf{32.3300}$\pm$14.4048&\textbf{42.4200}$\pm$1.5433&99.9900$\pm$4.2170&253.1900$\pm$3.4482&$--$&$--$\\
				RSHSVM$_{R}$&\textbf{8.7800}$\pm$0.3994&141.6300$\pm$0.8447&159.8700$\pm$0.0483&\textbf{9.2800}$\pm$0.3293&92.2200$\pm$2.2330&$--$&$--$\\
				SQSSVM&24.9600$\pm$0.3406&96.3200$\pm$7.3837&119.6200$\pm$1.0152&47.9200$\pm$0.3584&107.5600$\pm$1.2545&160.2400$\pm$1.5429&$--$\\
				QSMPM&$-$&$-$&$-$&$-$&$-$&$-$&$-$\\
				QSSVM$_{0/1}$&44.9700$\pm$12.9220&117.0700$\pm$5.5148&117.2800$\pm$5.8260&108.5500$\pm$13.9954&\textbf{55.4100}$\pm$15.8586&\textbf{119.0000}$\pm$0.0000&227.5800$\pm$4.7719\\
				\bottomrule
		\end{tabular}}
	\end{table*}

	\begin{figure}[!htbp]
		\centering 
		\includegraphics[ scale=0.4]{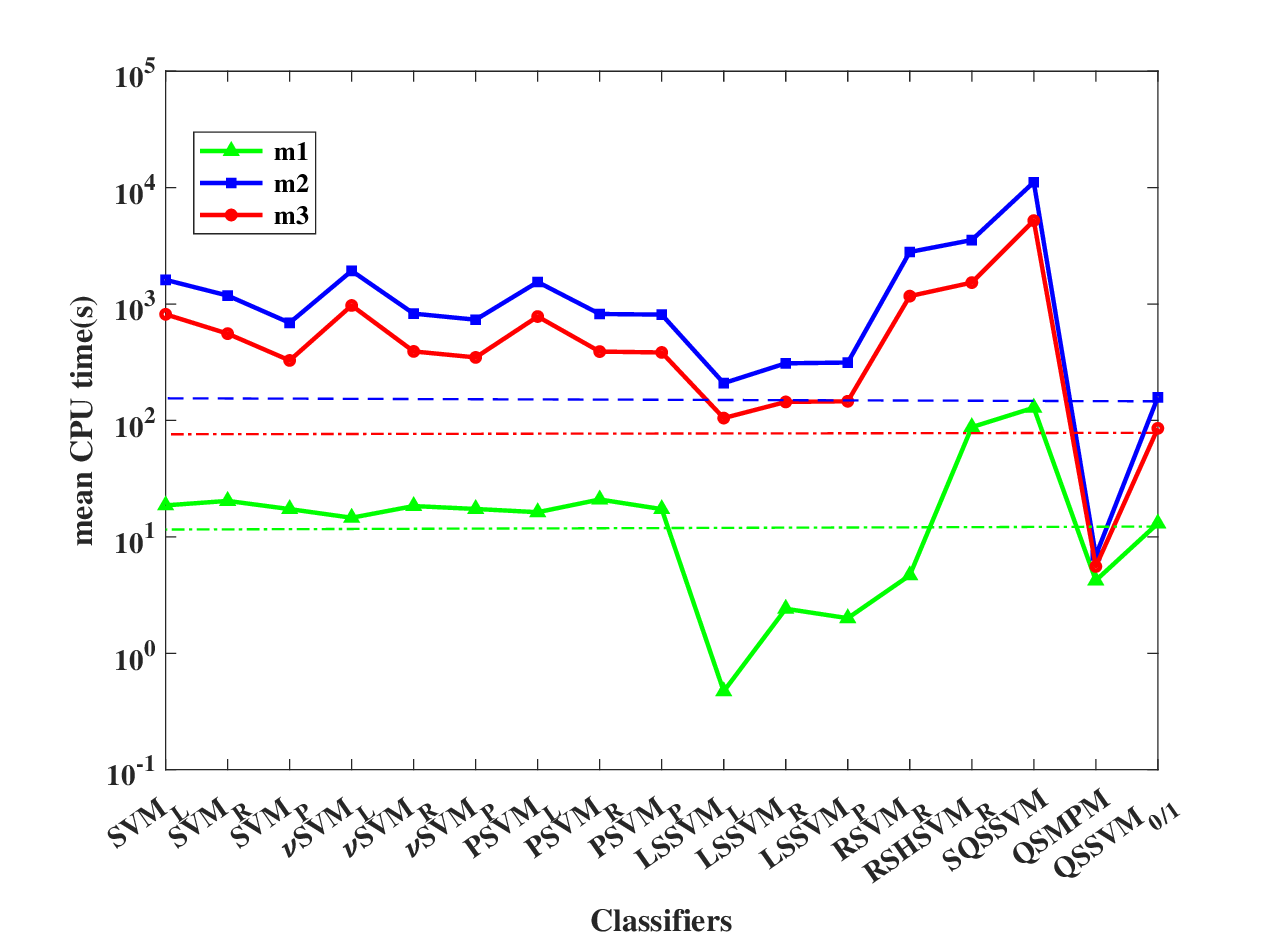}
		\centering
		\caption{Three mean values of CPU time}
		\label{figure_6}
	\end{figure}
		\Cref{figure_6} shows a line graph of the mean CPU time for each classifier. The labels ``m1'', ``m2'', and ``m3'' represent the mean CPU time on the former 7 small-scale datasets, the latter 7 large-scale datasets, and all datasets, respectively. It is worth noting that these three lines exhibit similar trends. Specifically, on the 7 small-scale datasets, the mean CPU time spent by our QSSVM $_{0/1}$ is higher than that of the 5 methods, but lower than that of the other methods. On the 7 large-scale datasets, the mean CPU time spent by our QSSVM $_{0/1}$ is higher than that of QSMPM and lower than that of the other methods. Across all datasets, our QSSVM$_{0/1}$ spends a lower mean CPU time than that of the 15 methods. Overall, our QSSVM$_{0/1}$ achieves superior computational efficiency with less mean CPU time compared to most other methods.
	
	In summary, the numerical experiments on 14 benchmark datasets indicate that our QSSVM$_{0/1}$ not only achieves strong sample sparsity, but also obtains high classification accuracy. In addition, the mean CPU time costs of our QSSVM$_{0/1}$ are relatively low. Therefore, our QSSVM$_{0/1}$ is feasible and effective.
	
		\subsection{Parameter analysis}
	
\begin{figure*}[!htbp]
		\centering
		\subfigure[Wireless]{
			\begin{minipage}[t]{0.2\linewidth}
				\centering
				\includegraphics[width=1.5in]{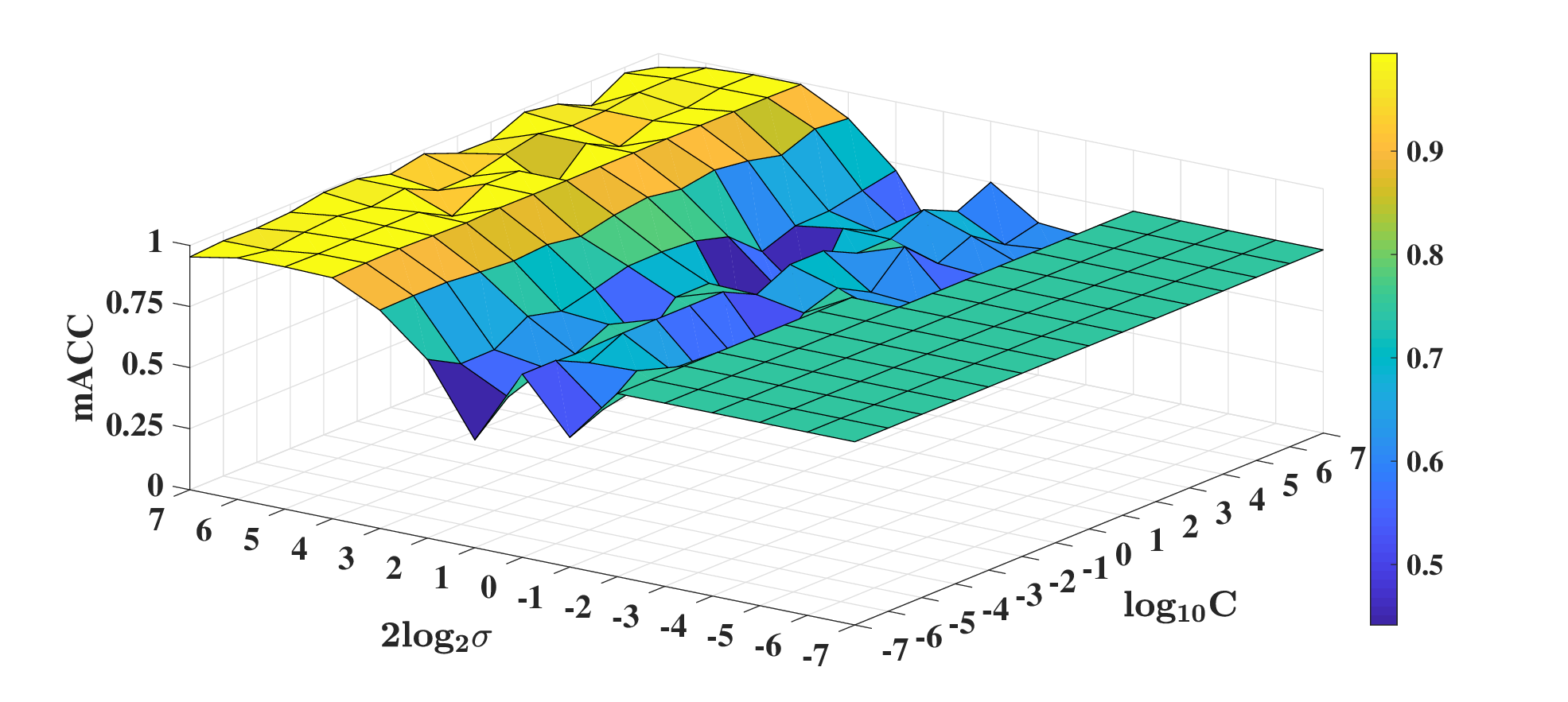}
			\end{minipage}
		}%
		\subfigure[Australian]{
			\begin{minipage}[t]{0.2\linewidth}
				\centering
				\includegraphics[width=1.5in]{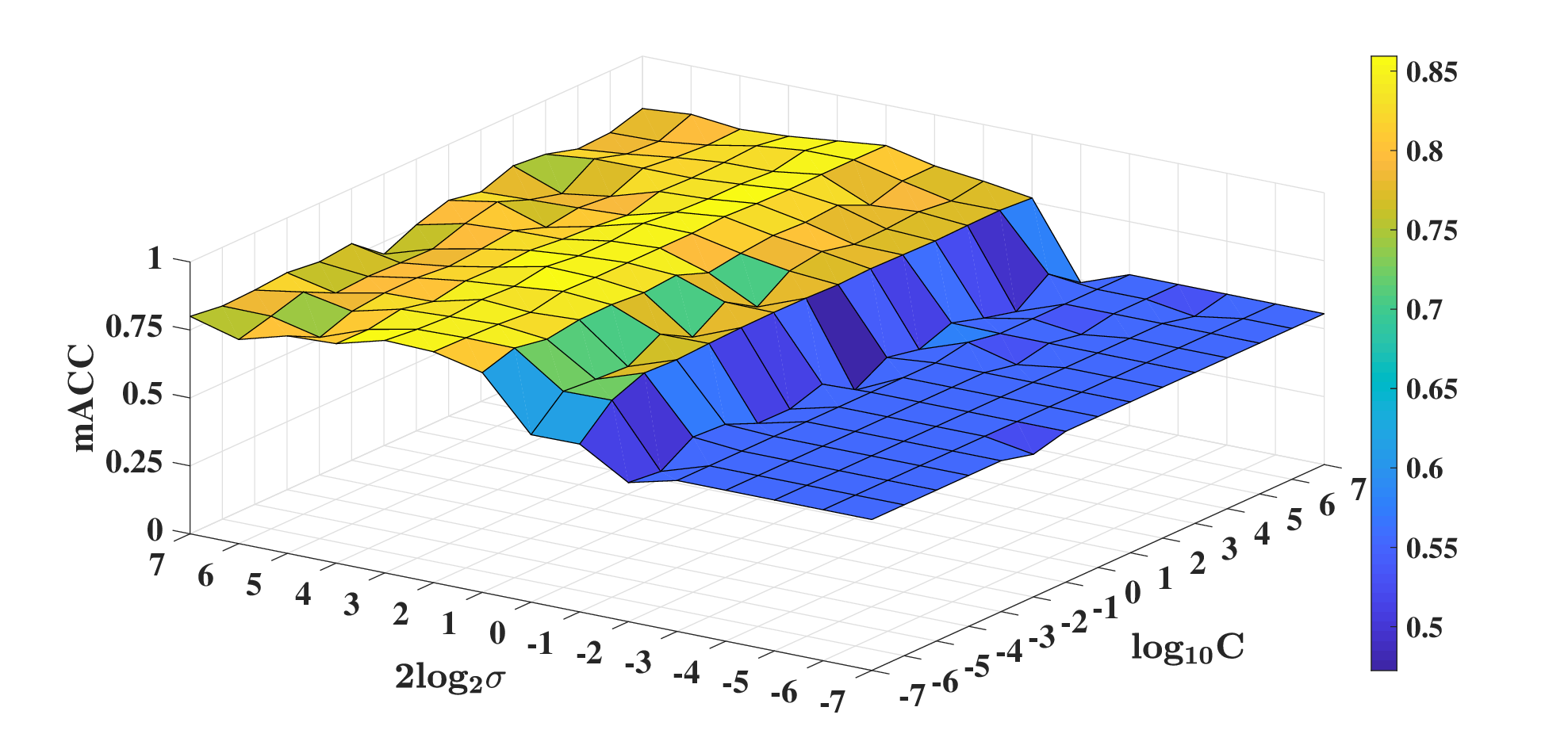}
			\end{minipage}
		}
		\subfigure[Banknote]{
			\begin{minipage}[t]{0.2\linewidth}
				\centering
				\includegraphics[width=1.6in]{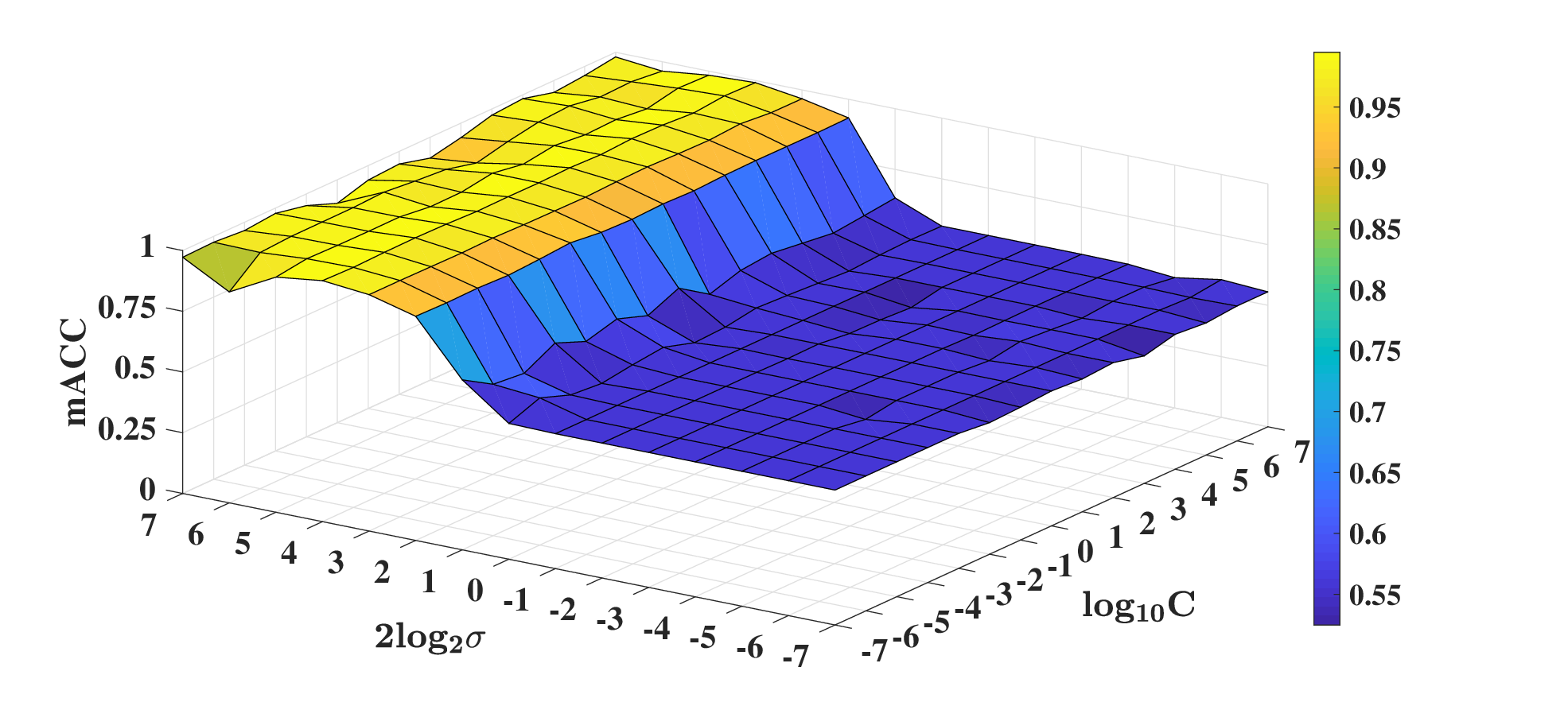}
			\end{minipage}
		}%
		\subfigure[Winequality]{
			\begin{minipage}[t]{0.2\linewidth}
				\centering
				\includegraphics[width=1.5in]{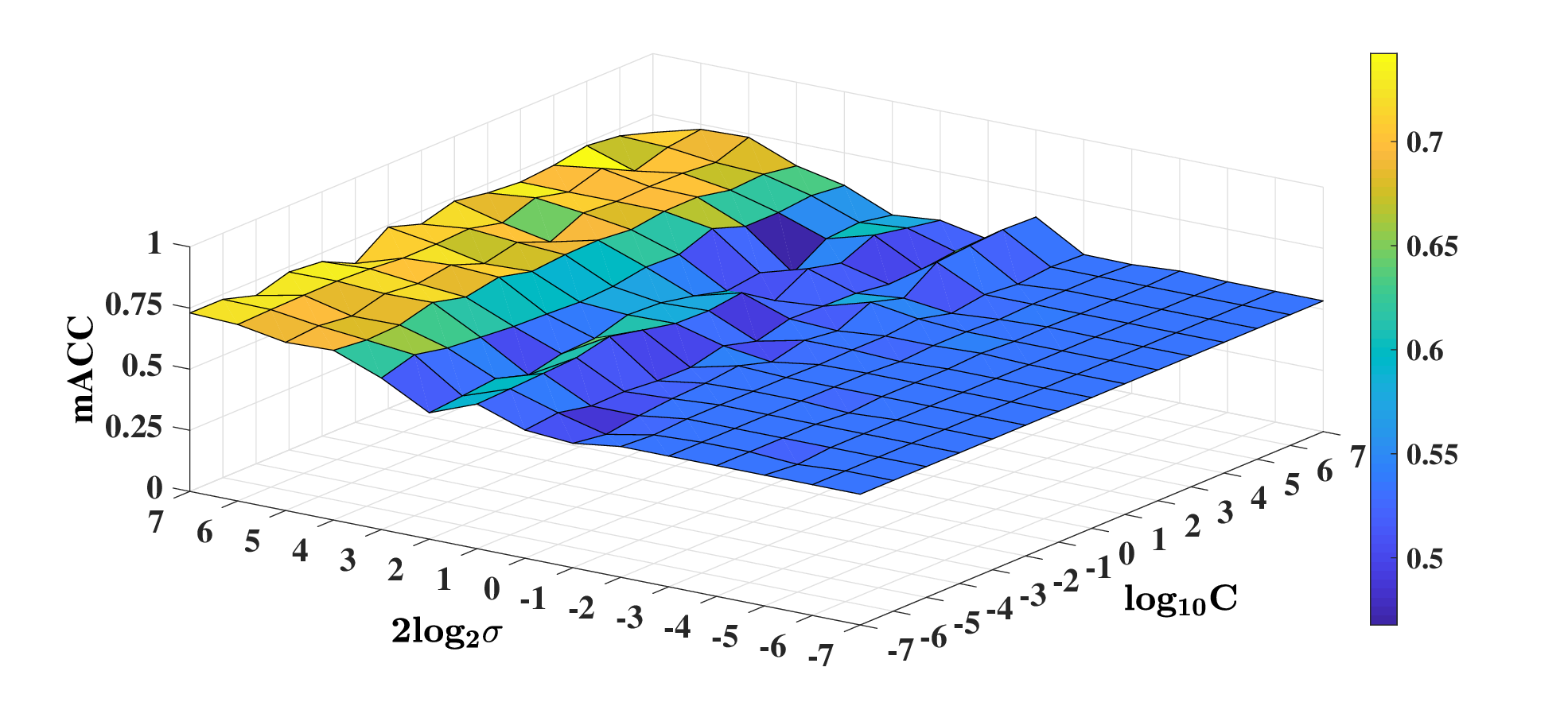}
			\end{minipage}
		}%
		\centering
		\caption{The mACC of QSSVM$_{0/1}$ versus the parameters $C$ and $\sigma$ on 4 benchmark datasets}\label{fcanshu}
	\end{figure*}

	To explore the effects of parameters $C$ and $\sigma$ on the mean accuracy of our QSSVM$_{0/1}$, grid plots for 4 benchmark datasets are shown in \Cref{fcanshu}, including Wireless, Australian, Banknote, Winequality datasets. \Cref{fcanshu} illustrates that the value of mACC tends to increase with increasing parameter $\sigma$ when parameter $C$ is fixed, while mACC is less sensitive to changes in the parameter $C$ when parameter $\sigma$ is fixed. These findings suggest that there is no need to spend too much time adjusting parameter $C$, and parameter $\sigma$ usually takes a larger number.
	
		\subsection{Statistical analysis}
\begin{figure}[!htbp]
		\subfigure[mACC]{
			\begin{minipage}[t]{0.5\linewidth}
				\centering
				\includegraphics[width=1.8in]{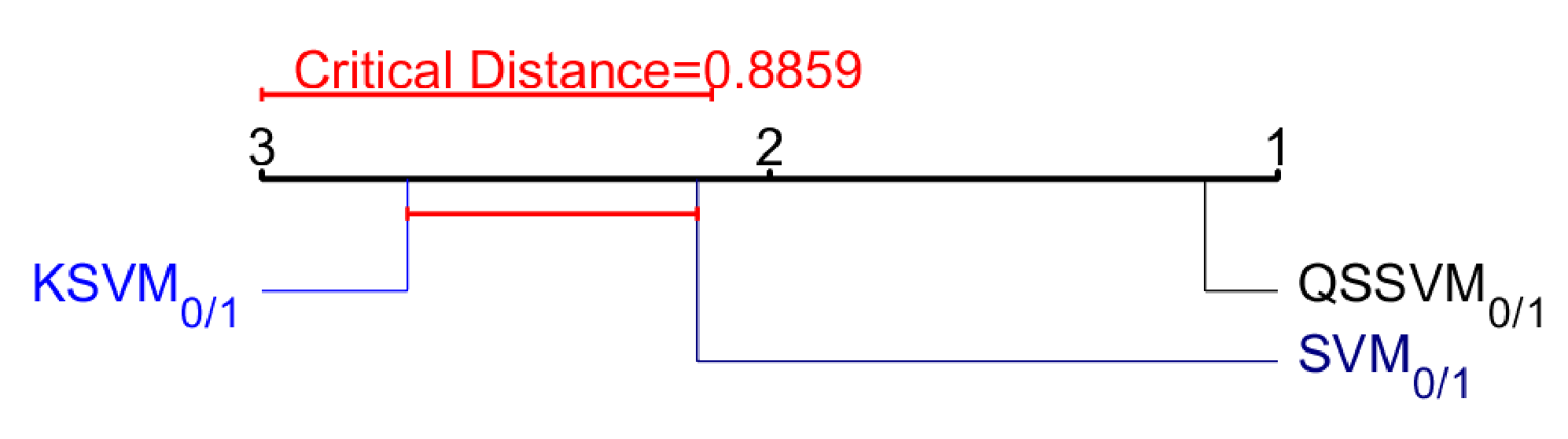}
			\end{minipage}
		}%
		\subfigure[mNSV]{
			\begin{minipage}[t]{0.5\linewidth}
				\centering
				\includegraphics[width=1.8in]{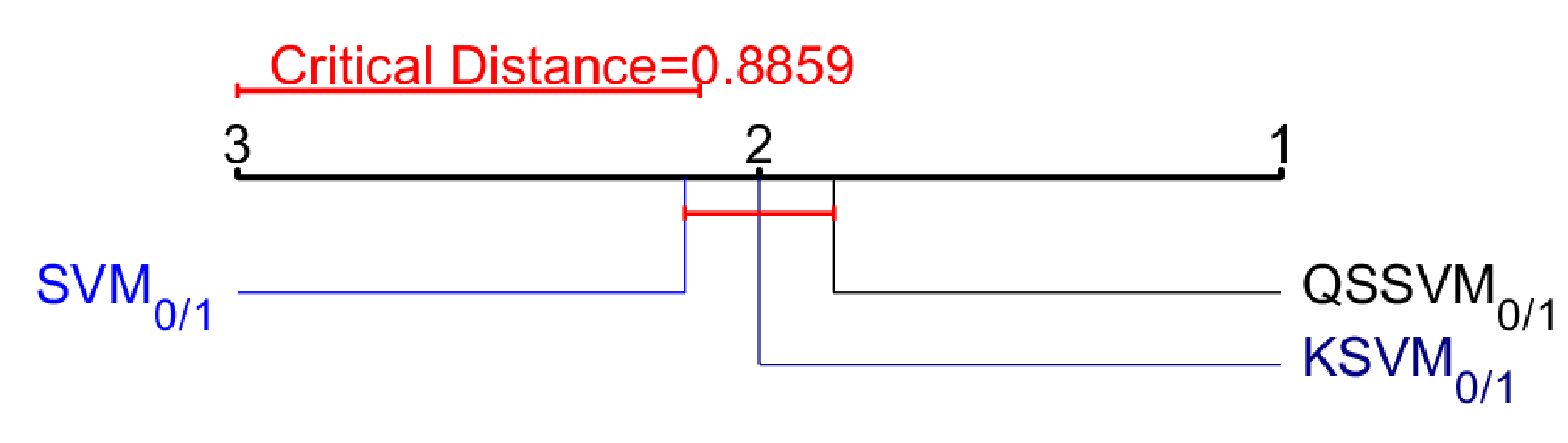}
			\end{minipage}
		}
		\centering
		\caption{The Friedman test and Nemenyi post-hoc test of our QSSVM$_{0/1}$ and SVM$_{0/1}$, KSVM$_{0/1}$ }\label{f7}
	\end{figure}

\begin{figure}[!htbp]
	\subfigure[mACC]{
		\begin{minipage}[t]{0.5\linewidth}
			\centering
			\includegraphics[trim=0 30 0 40, clip,width=1.8in]{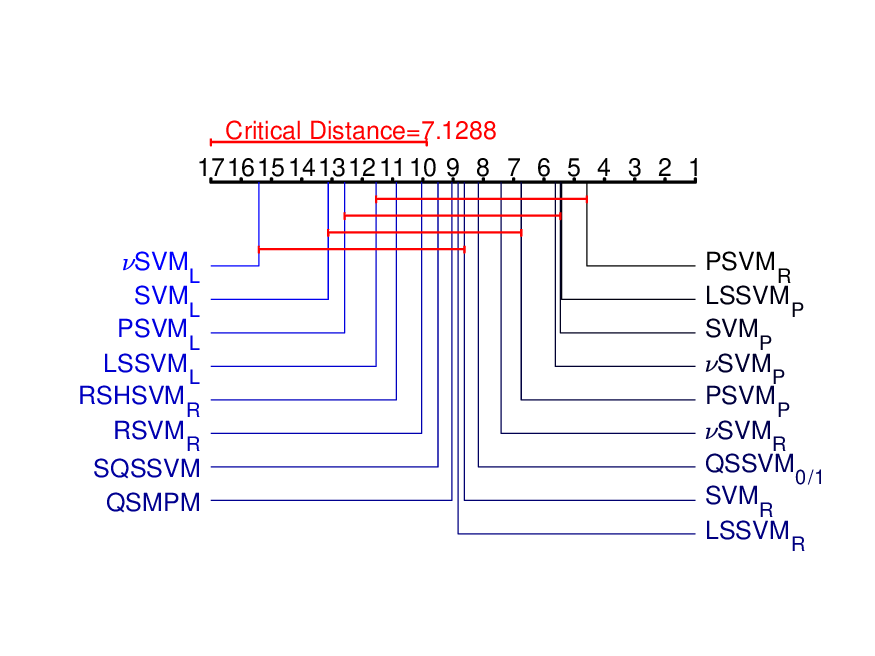}
		\end{minipage}
	}%
	\subfigure[mNSV]{
		\begin{minipage}[t]{0.5\linewidth}
			\centering
			\includegraphics[trim=0 30 0 40, clip,width=1.8in]{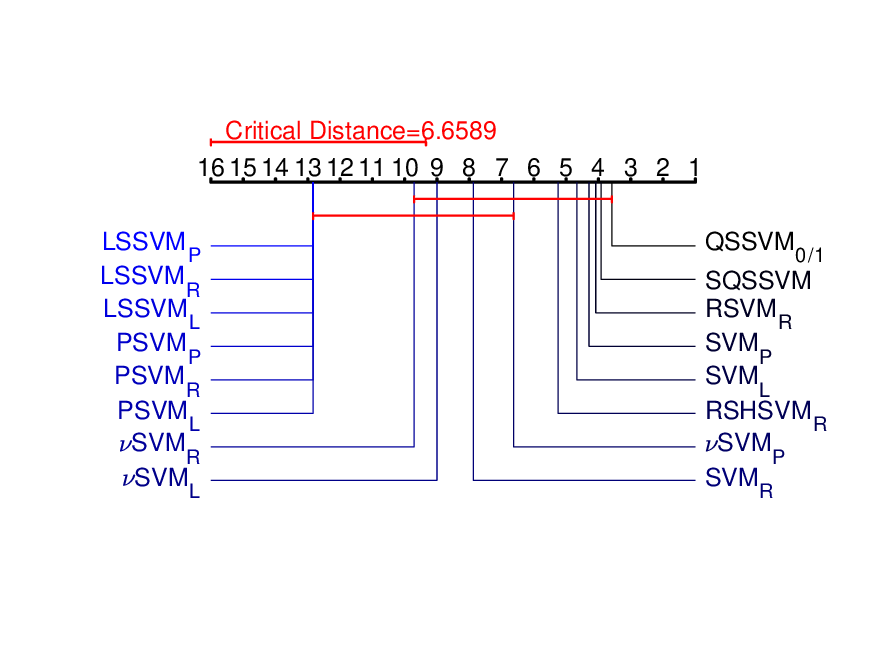}
		\end{minipage}
	}
	\centering
	\caption{The Friedman test and Nemenyi post-hoc test of our QSSVM$_{0/1}$ and other 16 methods}\label{f8}
\end{figure}

	In this subsection, Friedman test and post-hoc test are used to compare the performance of our QSSVM$_{0/1}$ with other methods. The results of the Friedman test and post-hoc test for the two evaluation criteria, including mACC and mNSV, are presented in \Cref{f7} and \Cref{f8}, respectively.
	
	For different methods, the Friedman test judges whether the null hypothesis of no significant difference at the significance level $\alpha=0.05$ is rejected. Then the post-hoc test is used to find out which methods differed significantly. To be specific, the Nemenyi test is used where the performance of two methods is significantly different if their average ranks over all datasets are larger than critical difference (CD), and the CD can be calculated by the following formula
	\begin{equation}
		CD=q_{\alpha}\sqrt{\frac{l(l+1)}{6h}},\label{CD}
	\end{equation}
	where $l$ is the number of methods and $h$ is the number of datasets.
	
	\Cref{f7} illustrates the Friedman test and Nemenyi post-hoc test on 14 benchmark datasets for our QSSVM$_{0/1}$ and 2 methods with 0-1 loss function, which are  SVM$_{0/1}$ and KSVM$_{0/1}$. The Friedman test in terms of mACC rejects the null hypothesis, but does not reject the null hypothesis in terms of mNSV. This indicates that these methods have significant difference in mACC, while there is no significant difference in mNSV. For $\alpha=0.05$, $q_{\alpha}=2.3440$. Thus the formula (\ref{CD}) yields $CD=0.8859$. Where the average ranks of each method are marked along the axis. The axis is turned so that the lowest (best) ranks are to the right. If there is no significant difference in the groups of methods they are connected by a red line. Statistically, our QSSVM$_{0/1}$ has more advantages in terms of mACC than that of SVM$_{0/1}$ and KSVM$_{0/1}$. For the evaluation criterion mNSV, QSSVM$_{0/1}$ is not significantly different from SVM$_{0/1}$ and KSVM$_{0/1}$.

The Friedman test and Nemenyi post-hoc test of our QSSVM$_{0/1}$ and 16 other methods on 12 benchmark datasets (excluding Waveform and Twonorm) are shown in \Cref{f8}. The Friedman test of the two evaluation criteria, mACC and mNSV, rejects the null hypothesis. Next, the Nemenyi post-hoc test is performed, for $\alpha=0.05$, we know in terms of the mACC, $q_{\alpha}=3.4580$, then the value of $CD=7.1288$ calculated by formula (\ref{CD}). And in terms of mNSV, $q_{\alpha}=3.4260$, the value of $CD=6.6589$. 
It is fairly easy to obtain that our QSSVM$_{0/1}$ ranks above average with respect to mACC and at the top with respect to mNSV. Therefore our QSSVM$_{0/1}$ has better effectiveness and sample sparsity.

	\section{Conclusion}\label{5}

In this paper, we develop a nonlinear kernel-free quadratic hyper-surface support vector machine with 0-1 loss function, namely QSSVM$_{0/1}$. It attempts to utilize a quadratic hyper-surface to divide the samples into two classes. By using the kernel-free trick, it reduces the computational time cost as it avoids choosing kernel functions and corresponding kernel parameters. And it has better interpretability compared to methods that use kernel functions. In addition, by introducing the 0-1 loss function to construct the optimization problem, our model achieves sample sparsity. Furthermore, according to the definitions of the P-stationary point and SVs, a new iterative algorithm based on the framework of the ADMM algorithm is designed to solve the optimization problem. It is worth noting that all SVs fall on the support hyper-surfaces. The relationship between the optimal solution to QSSVM$_{0/1}$ and P-stationary point is analyzed in the theoretical analysis. Numerical experiments indicate that our QSSVM$_{0/1}$ achieves higher classification accuracy, fewer number of SVs and less CPU time cost compared to most methods. 

In the future, to achieve feature selection, the $L_{1}$-norm $\|\cdot\|_{1}$ sparse regularization items will be added to the method to improve the performance of it. In addition, the 0-1 loss can also be applied to kernel-free soft quartic surface SVM model\cite{MR4188956} to achieve better classification performance.

\section*{Acknowledgments}
This work was supported by the National Natural Science Foundation of China (No.12061071), the Research Innovation Program for postgraduates of Xinjiang Uygur Autonomous Region (No.XJ2023G019), and the Xinjiang Key Laboratory of Applied Mathematics (No.XJDX1401).

\bibliographystyle{IEEEtran}
\bibliography{IEEEabrv,litrature}

	\appendix
\section*{Appendix A. Proofs}
\subsection*{A.1} \label{zhengming1}
\begin{proof}
	(1)Suppose that the matrix $\boldsymbol{D}$ defined by formula (\ref{D}) is full of column rank, there exists a generalized inverse matrix $\boldsymbol{D}^{+}$ of $\boldsymbol{D}$. Thus the constraint in optimization problem (\ref{yuan}) can be represented as $\boldsymbol{1}-\boldsymbol{u}=\boldsymbol{A}\widetilde{\boldsymbol{w}}+\boldsymbol{B}\boldsymbol{b}+c\boldsymbol{y}=\boldsymbol{D}(\widetilde{\boldsymbol{w}}; \boldsymbol{b}; c)$, so $[\widetilde{\boldsymbol{w}}; \boldsymbol{b}; c]=\boldsymbol{D}^{+}(\boldsymbol{1}-\boldsymbol{u})$. Based on this, the optimization problem (\ref{yuan}) can be transformed as a minimizer problem only containing single variable $\boldsymbol{u}$,
	\begin{equation}
		\min_{\boldsymbol{u}}\quad\dfrac{1}{2}\sum_{i=1}^{N}\|\boldsymbol{E}_{i}\boldsymbol{D}^{+}(\boldsymbol{u}-\boldsymbol{1})\|_{2}^
		{2}
		+C\|\boldsymbol{u}_{+}\|_{0},\label{singleu}
	\end{equation}
	where $\boldsymbol{D}, \boldsymbol{E}_{i},$ are defined by formulas (\ref{D}) and (\ref{Ei}). Suppose $\boldsymbol{u}^{*}$ is a globally optimal one to the optimization problem (\ref{yuan}), which is also a globally optimal solution to the optimization problem (\ref{singleu}) when the matrix $\boldsymbol{D}$ is full of column rank. Based on Lemma 3.1 in \cite{liu2024nonlinear}, denoting $L=\lambda_{\max}(\boldsymbol{H})$, and 
	\begin{equation}
		g(\boldsymbol{u}):=\dfrac{1}{2}\sum_{i=1}^{N}\|\boldsymbol{H}_{i}(\boldsymbol{u}-\boldsymbol{1})\|_{2}^
		{2},
	\end{equation}
	where $\boldsymbol{H}_{i}, \boldsymbol{H},$ are defined by the formulas (\ref{Hi}) and (\ref{H}) respectively, then we obtain $\nabla g(\boldsymbol{u})=\sum_{i=1}^{N}\boldsymbol{H}_{i}^\top\boldsymbol{H}_{i}(\boldsymbol{u}-\boldsymbol{1})$.
	Therefore, by the Lemma 3.1 in \cite{liu2024nonlinear} we have
	\begin{equation}
		\boldsymbol{u}^{*}=\operatorname{prox}_{\gamma C\|(\cdot)_{+}\|_{0}}(\boldsymbol{u}^{*}-\gamma\boldsymbol{H}(\boldsymbol{u}^{*}-\boldsymbol{1})),\label{uxing}
	\end{equation}
	for any $0<\gamma<1/\lambda_{\max}(\boldsymbol{H})$.
	Next, let
	\begin{equation}
		\boldsymbol{\lambda}^{*}:=\nabla g(\boldsymbol{u}^{*})=\boldsymbol{H}(\boldsymbol{u}^{*}-\boldsymbol{1}),\label{lamdaxing}
	\end{equation}
	then $\operatorname{prox}_{\gamma C\|(\cdot)_{+}\|_{0}}\left(\boldsymbol{u}^{*}-\gamma\boldsymbol{\lambda}^{*}\right)=\boldsymbol{u}^{*}$ is proved.
	
	Next, it follows from formula (\ref{lamdaxing}) that
	\begin{equation}\begin{aligned}
			-\boldsymbol{\lambda}^{*}
			&=\boldsymbol{H}(\boldsymbol{1}-\boldsymbol{u}^{*})\\
			&=\sum_{i=1}^{N}\boldsymbol{H}_{i}^\top\boldsymbol{H}_{i}(\boldsymbol{1}-\boldsymbol{u}^{*})\\
			&=\sum_{i=1}^{N}\boldsymbol{H}_{i}^\top\boldsymbol{E}_{i}\boldsymbol{D}^{+}(\boldsymbol{1}-\boldsymbol{u}^{*})\\
			&	=\sum_{i=1}^{N}\boldsymbol{H}_{i}^\top\boldsymbol{E}_{i}
			[\widetilde{\boldsymbol{w}}^{*};
			\boldsymbol{b}^{*};
			c^{*}].
		\end{aligned}
	\end{equation}
	It suffices to
	\begin{equation}
		\begin{aligned}
			-\boldsymbol{D}^\top\boldsymbol{\lambda}^{*}
			&=\sum_{i=1}^{N}\boldsymbol{D}^\top\boldsymbol{H}_{i}^\top\boldsymbol{E}_{i}
			[\widetilde{\boldsymbol{w}}^{*};
			\boldsymbol{b}^{*};
			c^{*}]\\
			&=\sum_{i=1}^{N}\boldsymbol{D}^\top\boldsymbol{D}^{+\top}\boldsymbol{E}_{i}^\top\boldsymbol{E}_{i}[\widetilde{\boldsymbol{w}}^{*};
			\boldsymbol{b}^{*};
			c^{*}]\\
			&=\sum_{i=1}^{N}\boldsymbol{E}_{i}^\top\boldsymbol{E}_{i}[\widetilde{\boldsymbol{w}}^{*};
			\boldsymbol{b}^{*};
			c^{*}]\\
			&=\left[
			\begin{array}{c}
				\sum_{i=1}^{N}\boldsymbol{M}_{i}^\top(\boldsymbol{M}_{i}\widetilde{\boldsymbol{w}^{*}} +\boldsymbol{b}^{*})\\
				\sum_{i=1}^{N}(\boldsymbol{M}_{i}\widetilde{\boldsymbol{w}}^{*} +\boldsymbol{b}^{*})\\
				0\end{array}\right],
		\end{aligned}
	\end{equation}
	where $\boldsymbol{D}^\top\boldsymbol{D}^{+}=\boldsymbol{I},   \boldsymbol{D}=[\boldsymbol{A}\, \boldsymbol{B}\, \boldsymbol{y}]\in\mathbb{R}^{N\times (\frac{n^{2}+n}{2}+n+1)}$.
	By the above formulas, the following formulas can be obtained
	$$
	\begin{cases}
		\begin{aligned}
			\sum_{i=1}^{N}\boldsymbol{M}_{i}^\top(\boldsymbol{M}_{i}\widetilde{\boldsymbol{w}}^{*} +\boldsymbol{b}^{*})+\boldsymbol{A}^\top\boldsymbol{\lambda}^{*}&=&\boldsymbol{0},\\
			\sum_{i=1}^{N}(\boldsymbol{M}_{i}\widetilde{\boldsymbol{w}}^{*} +\boldsymbol{b}^{*})+\boldsymbol{B}^\top\boldsymbol{\lambda}^{*}&=&\boldsymbol{0},\\
			\boldsymbol{y}^{\top}\boldsymbol{\lambda}^{*}&=&0.
		\end{aligned}
	\end{cases}
	$$
	
	Finally, the above conditions, the formula (\ref{uxing}) and the feasibility of $(\widetilde{\boldsymbol{w}}^{*}; \boldsymbol{b}^{*}; c^{*}; \boldsymbol{u}^{*})$ lead to formulas $(\ref{l1})$-$(\ref{l5})$. Therefore, the result (1) of Theorem \ref{theorem} is proved.
	
	(2) Suppose that $\Phi^{*}:=(\widetilde{\boldsymbol{w}}^{*}; \boldsymbol{b}^{*}; c^{*}; \boldsymbol{u}^{*})$ is a P-stationary point of the optimization problem (\ref{yuan}) and $\gamma>0$, then there is a $\boldsymbol{\lambda}^{*}\in\mathbb{R}^{N}$ such that $(\Phi^{*}; \boldsymbol{\lambda}^{*})$ satisfies formulas (\ref{l1})-(\ref{l5}).
	
	Let $\Theta$ be the feasible region of the optimization problem (\ref{yuan}), namely
	\begin{equation}
		\Theta:=\{\Phi:\boldsymbol{u}+\boldsymbol{A}\widetilde{\boldsymbol{w}}+\boldsymbol{B}\boldsymbol{b}+c\boldsymbol{y}=\boldsymbol{1}\},
	\end{equation}
	Moreover, $\|\boldsymbol{u}_{+}\|_{0}$ is lower semi-continuous at $\Phi^{*}\in\Theta$, then by the Proposition 4.3 of \cite{Mordukhovich}, there is a neighborhood $U(\Phi^{*}, \delta_{1})$ of $\Phi^{*}\in\Theta$ with $\delta_{1}>0$ such that
	\begin{equation}
		\|\boldsymbol{u}_{+}\|_{0}>\|\boldsymbol{u}^{*}_{+}\|_{0}-\dfrac{1}{2}, \quad\forall\Phi\in\Theta\cap U(\Phi^{*}, \delta_{1})\label{U}.
	\end{equation}
	Besides,  $ \sum_{i=1}^{N}\|\boldsymbol{M}_{i}\widetilde{\boldsymbol{w}} +\boldsymbol{b}\|_{2}^
	{2}$ is locally lipschitz continuous, there exists a neighborhood $U(\Phi^{*}, \delta_{2})$ of $\Phi^{*}\in\Theta$ with $\delta_{2}>0$ such that
	\begin{equation}
		\begin{aligned}
			&|\sum_{i=1}^{N}\|\boldsymbol{M}_{i}\widetilde{\boldsymbol{w}} +\boldsymbol{b}\|_{2}^
			{2}-\sum_{i=1}^{N}\|\boldsymbol{M}_{i}\widetilde{\boldsymbol{w}}^{*} +\boldsymbol{b}^{*}\|_{2}^
			{2}|\leq2C, \\
			& \quad \forall\Phi\in\Theta\cap U(\Phi^{*}, \delta_{2}).\label{z'gz}
		\end{aligned}
	\end{equation}
	
	Denote $\delta:=\min\{\delta_{1}, \delta_{2}\}$, next, we prove that $\Phi^{*}$ is a local minimizer for the optimization problem (\ref{yuan}). Namely, there exists a neighborhood $U(\Phi^{*}, \delta)$ of $\Phi^{*}\in\Theta$ with $\delta>0$ such that
	\begin{equation}
		\begin{aligned}
			&	\dfrac{1}{2}\sum_{i=1}^{N}\|\boldsymbol{M}_{i}\widetilde{\boldsymbol{w}}^{*} +\boldsymbol{b}^{*}\|_{2}^
			{2}+C\|\boldsymbol{u}^{*}_{+}\|_{0}\\
			\leq&\dfrac{1}{2}\sum_{i=1}^{N}\|\boldsymbol{M}_{i}\widetilde{\boldsymbol{w}} +\boldsymbol{b}\|_{2}^
			{2}
			+C\|\boldsymbol{u}_{+}\|_{0}, \\
			&	\quad\forall\Phi\in\Theta\cap U(\Phi^{*}, \delta).\label{jubujizhi}
		\end{aligned}
	\end{equation}
	For this purpose, let $\Gamma_{*}:=\{i: u_{i}^{*}=0\}, \overline{\Gamma}_{*}:=\mathbb{N}_{N}\setminus\Gamma_{*}$. According to formulas (\ref{l5}) and (\ref{popetor1}) we have
	\begin{equation}
		\begin{aligned}
			-\sqrt{2C/\gamma}\leq\lambda_{i}^{*}&\leq0, u_{i}^{*}=0, \quad\forall i\in\Gamma_{*},\\
			\lambda_{i}^{*}&=0, u_{i}^{*}\neq0, \quad\forall i\in\overline{\Gamma}_{*}.\label{guanxi1}
		\end{aligned}
	\end{equation}
	Based on these, a local region $\Theta_{1}$ of $\Theta$ is  considered, i.e.,
	\begin{equation}
		\Theta_{1}:=\Theta\cap\{\Phi: u_{i}\leq0, i\in\Gamma_{*}\}.
	\end{equation}
	The proof of inequality (\ref{jubujizhi}) is divided into the following two cases:
	
	(i) $\Phi\in\Theta_{1}\subseteq\Theta$ and $\Phi\in U(\Phi^{*}, \delta)$. It is to get that $\Phi^{*}\in\Theta_{1}$ by formulas (\ref{l1})-(\ref{l5}). Then for any $\Phi\in\Theta_{1}$, we have
	\begin{eqnarray}
		u_{i}\leq0, i\in\Gamma_{*}, \label{ui}\\
		\boldsymbol{u}+\boldsymbol{A}\widetilde{\boldsymbol{w}}+\boldsymbol{B}\boldsymbol{b}+c\boldsymbol{y}=\boldsymbol{1},
	\end{eqnarray}
	which and formula (\ref{l4}) suffice to
	\begin{equation}
		-\boldsymbol{A}(\widetilde{\boldsymbol{w}}-\widetilde{\boldsymbol{w}}^{*})=(\boldsymbol{u}-\boldsymbol{u}^{*})+\boldsymbol{B}(\boldsymbol{b}-\boldsymbol{b}^{*})+(c-c^{*})\boldsymbol{y}.\label{AZUCY}
	\end{equation}
	The following chain of inequalities hold for any $\Phi\in\Theta_{1}$,
	\begin{equation}
		\begin{array}{cl}
			&\sum_{i=1}^{N}\|\boldsymbol{M}_{i}\widetilde{\boldsymbol{w}} +\boldsymbol{b}\|_{2}^
			{2}-\sum_{i=1}^{N}\|\boldsymbol{M}_{i}\widetilde{\boldsymbol{w}}^{*} +\boldsymbol{b}^{*}\|_{2}^
			{2}\\
			=&\sum_{i=1}^{N}\|\boldsymbol{M}_{i}\widetilde{\boldsymbol{w}} +\boldsymbol{b}-\boldsymbol{M}_{i}\widetilde{\boldsymbol{w}}^{*} -\boldsymbol{b}^{*}+\boldsymbol{M}_{i}\widetilde{\boldsymbol{w}}^{*} +\boldsymbol{b}^{*}\|_{2}^
			{2}\\
			&-\sum_{i=1}^{N}\|\boldsymbol{M}_{i}\widetilde{\boldsymbol{w}}^{*} +\boldsymbol{b}^{*}\|_{2}^
			{2}\\
			=&\sum_{i=1}^{N}\|\boldsymbol{M}_{i}\widetilde{\boldsymbol{w}} +\boldsymbol{b}-\boldsymbol{M}_{i}\widetilde{\boldsymbol{w}}^{*} -\boldsymbol{b}^{*}\|_{2}^
			{2}\\
			&+2\sum_{i=1}^{N}\langle\boldsymbol{M}_{i}\widetilde{\boldsymbol{w}} +\boldsymbol{b}-\boldsymbol{M}_{i}\widetilde{\boldsymbol{w}}^{*} -\boldsymbol{b}^{*}, \boldsymbol{M}_{i}\widetilde{\boldsymbol{w}}^{*}+\boldsymbol{b}^{*}\rangle\\
			\geq&2\sum_{i=1}^{N}\langle\boldsymbol{M}_{i}\widetilde{\boldsymbol{w}} +\boldsymbol{b}-\boldsymbol{M}_{i}\widetilde{\boldsymbol{w}}^{*} -\boldsymbol{b}^{*}, \boldsymbol{M}_{i}\widetilde{\boldsymbol{w}}^{*}+\boldsymbol{b}^{*}\rangle\\
			\stackrel{(\ref{l1}), (\ref{l2})}{=}&-2(\widetilde{\boldsymbol{w}}-\widetilde{\boldsymbol{w}}^{*})^\top\boldsymbol{A}^{\top}\boldsymbol{\lambda}^{*}-2\boldsymbol{\lambda}^{\top*}\boldsymbol{B}(\boldsymbol{b}-\boldsymbol{b}^{*}) \\
			\stackrel{(\ref{AZUCY})}{=}&2(\boldsymbol{u}-\boldsymbol{u}^{*})^{\top} \boldsymbol{\lambda}^{*}+
			2(c-c^{*})\boldsymbol{y}^{\top} \boldsymbol{\lambda}^{*}\\
			\stackrel{(\ref{l3})}{=}&2(\boldsymbol{u}-\boldsymbol{u}^{*})^{\top} \boldsymbol{\lambda}^{*}\\
			=&2(\boldsymbol{u}_{\Gamma_{*}}-\boldsymbol{u}_{\Gamma_{*}}^{*})^{\top} \boldsymbol{\lambda}_{\Gamma_{*}}^{*}+
			2(\boldsymbol{u}_{\overline{\Gamma}_{*}}-\boldsymbol{u}_{\overline{\Gamma}_{*}}^{*})^{\top} \boldsymbol{\lambda}_{\overline{\Gamma}_{*}}^{*}\\
			\stackrel{(\ref{guanxi1})}{=}&2\boldsymbol{u}_{{\Gamma}_{*}}^{\top} \boldsymbol{\lambda}_{\Gamma_{*}}^{*}\\
			\stackrel{(\ref{guanxi1}), (\ref{ui})}{\geq}&0.\label{1111}
		\end{array}
	\end{equation}
	Since $\|u_{+}\|_{0}$ can only take values from $\{0, 1, \ldots, N\}$, this together with formula (\ref{U}) allows us to conclude that
	\begin{equation}
		\|\boldsymbol{u}_{+}\|_{0}\geq\|\boldsymbol{u}^{*}_{+}\|_{0}, \quad\forall\Phi\in\Theta\cap U(\Phi^{*}, \delta_{1}).\label{2222}
	\end{equation}
	Hence, for any $\Phi\in\Theta_{1}\cap U(\Phi^{*}, \delta)\subseteq\Theta_{1}\cap U(\Phi^{*}, \delta_{1})$, then formulas (\ref{1111}) and (\ref{2222}) lead to
	\begin{equation}
		\dfrac{1}{2}\|\boldsymbol{M}_{i}\widetilde{\boldsymbol{w}}^{*} +\boldsymbol{b}^{*}\|_{2}^
		{2}
		+C\|\boldsymbol{u}^{*}_{+}\|_{0}\leq\dfrac{1}{2}\|\boldsymbol{M}_{i}\widetilde{\boldsymbol{w}} +\boldsymbol{b}\|_{2}^
		{2}
		+C\|\boldsymbol{u}_{+}\|_{0}.\label{guanjian1}
	\end{equation}
	
	(ii) $\Phi\in(\Theta\setminus\Theta_{1})$ and $\Phi\in U(\Phi^{*}, \delta)$. On the basis of $\Phi\in(\Theta\setminus\Theta_{1})$, there exists $i_{0}\in\Gamma_{*}$ with $u_{i0}^{*}=0$, but $u_{i0}>0$, which indicates that $\|(u_{i_0}^*)_{+}\|_0=0, \|(u_{i_0})_{+}\|_0=1$. By $\Phi\in U(\Phi^{*}, \delta)$ and formula (\ref{2222}) we have
	\begin{equation}
		\|\boldsymbol{u}_{+}\|_{0}\geq\|\boldsymbol{u}^{*}_{+}\|_{0}+1.
	\end{equation}
	This together with formula (\ref{z'gz}) obtains that for any $\Phi\in(\Theta\setminus\Theta_{1})\cap U(\Phi^{*}, \delta)$,
	\begin{equation}
		\begin{aligned}
			&\dfrac{1}{2}\|\boldsymbol{M}_{i}\widetilde{\boldsymbol{w}}^{*} +\boldsymbol{b}^{*}\|_{2}^
			{2}+C\|\boldsymbol{u}^{*}_{+}\|_{0}\\
			\leq&\dfrac{1}{2}\|\boldsymbol{M}_{i}\widetilde{\boldsymbol{w}}^{*} +\boldsymbol{b}^{*}\|_{2}^
			{2}+C\|\boldsymbol{u}_{+}\|_{0}-C\\
			\leq&\dfrac{1}{2}\|\boldsymbol{M}_{i}\widetilde{\boldsymbol{w}} +\boldsymbol{b}\|_{2}^
			{2}+C\|\boldsymbol{u}_{+}\|_{0}.\label{guanjian2}
		\end{aligned}
	\end{equation}
	Combining formulas (\ref{guanjian1}) and (\ref{guanjian2}), we obtain that $\Phi^{*}$ is a local minimizer of the optimization problem (\ref{yuan}) in a local region $\Theta\cap U(\Phi^{*}, \delta)$. Thus the conclusion is proved.
\end{proof}
\subsection*{A.2} \label{zhengming2}
\begin{proof}
	Since $T_{k}\in\mathbb{N}_{N}$, for sufficient large $k$, there is a subset $J\subseteq\{1, 2, 3, \ldots\}$ such that
	\begin{equation}
		T_{j}\equiv:T, \quad\forall j\in J.
	\end{equation}
	For simplicity, denote the sequence $\Psi^{k}:=(\widetilde{\boldsymbol{w}}^{k}, \boldsymbol{b}^{k}, c^{k}, \boldsymbol{u}^{k}, \boldsymbol{\lambda}^{k})$ and its limit point $\Psi^{*}:=(\widetilde{\boldsymbol{w}}^{*}, \boldsymbol{b}^{*}, c^{*}, \boldsymbol{u}^{*}, \boldsymbol{\lambda}^{*})$, namely $\{\Psi^{k}\}\rightarrow\Psi^{*}$. This also indicates that $\{\Psi^{j}\}_{j\in J}\rightarrow\Psi^{*}$ and $\{\Psi^{j+1}\}_{j\in J}\rightarrow\Psi^{*}$. 
	
	Taking the limit along with $J$ of formula (\ref{lamda}), namely, $k\in J, k\rightarrow\infty$,  we have
	\begin{equation}
		\begin{cases}
			\boldsymbol{\lambda}_{T}^{*}=\boldsymbol{\lambda}_{T}^{*}+\eta\sigma\boldsymbol{\varpi}_{T}^{*},\\
			\boldsymbol{\lambda}_{\overline{T}}^{*}=\boldsymbol{0}, \label{l*}
		\end{cases}
	\end{equation}
	where $\boldsymbol{\varpi}:=\boldsymbol{u}-\boldsymbol{1}+\boldsymbol{A}\widetilde{\boldsymbol{w}}+\boldsymbol{B}\boldsymbol{b}
	+c\boldsymbol{y}$, which derives $\boldsymbol{\varpi}_{T}^{*}=\boldsymbol{0}$. Taking the limit along with $J$ of formula (\ref{updateuk+1}) and $\boldsymbol{v}^{k}$ respectively yields
	\begin{equation}
		\begin{aligned}
			\boldsymbol{v}^{*}&= \boldsymbol{1}-\boldsymbol{A}\widetilde{\boldsymbol{w}}^{*}-\boldsymbol{B}\boldsymbol{b}^{*}-c^{*}\boldsymbol{y}-\boldsymbol{\lambda}^{*}/{\sigma}\\
			&= \boldsymbol{1}-\boldsymbol{A}\widetilde{\boldsymbol{w}}^{*}-\boldsymbol{B}\boldsymbol{b}^{*}-c^{*}\boldsymbol{y}-\boldsymbol{u}^{*}+\boldsymbol{u}^{*}-\boldsymbol{\lambda}^{*}/{\sigma}\\
			&= -\boldsymbol{\varpi}_{T}^{*}+\boldsymbol{u}^{*}-\boldsymbol{\lambda}^{*}/{\sigma},\label{v*}
		\end{aligned}
	\end{equation}
	and thus,
	\begin{eqnarray}
		\boldsymbol{u}^{*}_{T}&=&\boldsymbol{0},\label{ut}\\
		\boldsymbol{u}^{*}_{\overline{T}}&=&\boldsymbol{v}^{*}_{\overline{T}},\label{ut1}\\
		&\stackrel{(\ref{v*})}{=}&-\boldsymbol{\varpi}_{\overline{T}}^{*}+\boldsymbol{u}_{\overline{T}}^{*}-\boldsymbol{\lambda}_{\overline{T}}^{*}/{\sigma},\\
		&\stackrel{(\ref{l*})}{=}&-\boldsymbol{\varpi}_{\overline{T}}^{*}+\boldsymbol{u}_{\overline{T}}^{*}.
	\end{eqnarray}
	These indicate that $\boldsymbol{\varpi}_{\overline{T}}^{*}=\boldsymbol{0}$, therefore, $\boldsymbol{\varpi}^{*}=\boldsymbol{0}$, where $\boldsymbol{\varpi}^{*}=\boldsymbol{u}^{*}-\boldsymbol{1}+\boldsymbol{A}\widetilde{\boldsymbol{w}}^{*}+\boldsymbol{B}\boldsymbol{b}^{*}
	+c^{*}\boldsymbol{y}$. Then we have 
	\begin{equation}	\boldsymbol{u}^{*}+\boldsymbol{A}\widetilde{\boldsymbol{w}}^{*}+\boldsymbol{B}\boldsymbol{b}^{*}+c^{*}\boldsymbol{y}-\boldsymbol{1}=\boldsymbol{0},\label{u1}
	\end{equation}
	and  $\boldsymbol{v}^{*}=\boldsymbol{u}^{*}-\boldsymbol{\lambda}^{*}/{\sigma}$ by formula (\ref{v*}), which together with formulas (\ref{ut}), (\ref{ut1}) and the proximal operator of 0-1 loss function \cite{15921057020221001} indicates
	\begin{equation}
		\begin{aligned}
			\boldsymbol{u}^{*}&=\operatorname{prox}_{\frac{C}{\sigma}\|(\cdot)_{+}\|_{0}}\left(\boldsymbol{v}^{*}\right)\\
			&=\operatorname{prox}_{\frac{C}{\sigma}\|(\cdot)_{+}\|_{0}}\left(\boldsymbol{u}^{*}-\boldsymbol{\lambda}^{*}/\sigma\right).\label{linjin}
		\end{aligned}
	\end{equation}
	
	Taking the limit along with $J$ of formula (\ref{updatew}) leads to
	\begin{equation}
		\begin{aligned}
			&(\boldsymbol{G}+\sigma[\boldsymbol{A}_{T}\:\boldsymbol{B}_{T}]^{\top}[\boldsymbol{A}_{T}\:\boldsymbol{B}_{T}])[\widetilde{\boldsymbol{w}}^{*}; \boldsymbol{b}^{*}]\\
			=&-\sigma [\boldsymbol{A}_{T}\:\boldsymbol{B}_{T}]^{\top}(\boldsymbol{u}^{*}_{T}+c^{*}\boldsymbol{y}-\boldsymbol{1}+\boldsymbol{\lambda}^{*}_{T}/\sigma)
			\\
			=&-\sigma[\boldsymbol{A}_{T}\:\boldsymbol{B}_{T}]^{\top}(\boldsymbol{\varpi}^{*}_{T}-[\boldsymbol{A}_{T}\:\boldsymbol{B}_{T}][\widetilde{\boldsymbol{w}}^{*}; \boldsymbol{b}^{*}]+\boldsymbol{\lambda}^{*}_{T}/\sigma)
			\\
			=&-\sigma[\boldsymbol{A}_{T}\:\boldsymbol{B}_{T}]^{\top}(-[\boldsymbol{A}_{T}\:\boldsymbol{B}_{T}][\widetilde{\boldsymbol{w}}^{*}; \boldsymbol{b}^{*}]+\boldsymbol{\lambda}^{*}_{T}/\sigma).\label{jixian}
		\end{aligned}
	\end{equation}
	Formula (\ref{jixian}) can be rewritten as
	$$\begin{aligned}
		&\boldsymbol{G}[\widetilde{\boldsymbol{w}}^{*};\boldsymbol{b}^{*}]+\sigma[\boldsymbol{A}_{T}\:\boldsymbol{B}_{T}]^{\top}[\boldsymbol{A}_{T}\:\boldsymbol{B}_{T}][\widetilde{\boldsymbol{w}}^{*}; \boldsymbol{b}^{*}]\\
		=&\sigma[\boldsymbol{A}_{T}\:\boldsymbol{B}_{T}]^{\top}[\boldsymbol{A}_{T}\:\boldsymbol{B}_{T}][\widetilde{\boldsymbol{w}}^{*}; \boldsymbol{b}^{*}]-[\boldsymbol{A}_{T}\:\boldsymbol{B}_{T}]^{\top}\boldsymbol{\lambda}_{T}^{*},\end{aligned}
	$$
	then, $$\boldsymbol{G}[\widetilde{\boldsymbol{w}}^{*}; \boldsymbol{b}^{*}]
	=-[\boldsymbol{A}_{T}\:\boldsymbol{B}_{T}]^{\top}\boldsymbol{\lambda}_{T}^{*}
	\stackrel{(\ref{l*})}{=}-[\boldsymbol{A}\:\boldsymbol{B}]^{\top}\boldsymbol{\lambda}^{*}
	,$$
	so, we have
	$\boldsymbol{G}[\widetilde{\boldsymbol{w}}^{*}; \boldsymbol{b}^{*}]+[\boldsymbol{A}_{T}\:\boldsymbol{B}_{T}]^{\top}\boldsymbol{\lambda}^{*}=\boldsymbol{0}$,
	namely,
	
	\begin{eqnarray}
		\sum_{i=1}^{N}\boldsymbol{M}_{i}^\top(\boldsymbol{M}_{i}\widetilde{\boldsymbol{w}}^{*} +\boldsymbol{b}^{*})+\boldsymbol{A}^\top\boldsymbol{\lambda}^{*}&=&\boldsymbol{0},\label{GW}\\
		\sum_{i=1}^{N}(\boldsymbol{M}_{i}\widetilde{\boldsymbol{w}}^{*} +\boldsymbol{b}^{*})+\boldsymbol{B}^\top\boldsymbol{\lambda}^{*}&=&\boldsymbol{0},\label{Ib}.
	\end{eqnarray}
	
Taking the limit along with $J$ of formula (\ref{updateck+1})
$$
\begin{aligned}
	c^{*}
	&=-\boldsymbol{y}^{\top} (\boldsymbol{A}\widetilde{\boldsymbol{w}}^{*}+\boldsymbol{B}\boldsymbol{b}^{*}-\boldsymbol{1}+\boldsymbol{u}^{*}+\boldsymbol{\lambda}^{*}/\sigma)/N\\
	&=-\boldsymbol{y}^{\top} (\boldsymbol{A}\widetilde{\boldsymbol{w}}^{*}+\boldsymbol{B}\boldsymbol{b}^{*}-\boldsymbol{1}+\boldsymbol{u}^{*}+c^{*}\boldsymbol{y}-c^{*}\boldsymbol{y}+\boldsymbol{\lambda}^{*}/\sigma)/N\\
	&=-\boldsymbol{y}^{\top} (\boldsymbol{\varpi}^{*}-c^{*}\boldsymbol{y}+\boldsymbol{\lambda}^{*}/\sigma)/N\\
	&=-\boldsymbol{y}^{\top} (-c^{*}\boldsymbol{y}+\boldsymbol{\lambda}^{*}/\sigma)/N\\
	&=(c^{*}\boldsymbol{y}^{\top}\boldsymbol{y}-\boldsymbol{y}^{\top}\boldsymbol{\lambda}^{*}/{\sigma})/N\\
	&=c^{*}-\boldsymbol{y}^{\top}\boldsymbol{\lambda}^{*}/(\sigma N),
\end{aligned}
$$
which indicates that \begin{equation}\boldsymbol{y}^{\top}\boldsymbol{\lambda}^{*}=0.\label{yl}
\end{equation}
Based on the results of equations (\ref{u1}), (\ref{linjin}), (\ref{GW}), (\ref{Ib}) and (\ref{yl}), we can conclude that $(\widetilde{\boldsymbol{w}}^{*}; \boldsymbol{b}^{*}; c^{*}; \boldsymbol{u}^{*})$ is a P-stationary of the optimization problem (\ref{yuan}), where $\gamma=1/\sigma$. Then by Theorem \ref{theorem} (2), it is a locally optimal solution to the optimization problem (\ref{yuan}). Therefore, the proof is complete.
\end{proof}

\end{document}